\input amstex
\documentstyle{amsppt}
%
% Author:   Ruslan Abdulovich Sharipov
% Title:    Dynamical systems admitting the normal shift
% Comment:  Thesis for the degree of Doctor of Sciences in Russia
% Address:  Rabochaya street 5, Ufa, 450003, Russia
% E-mail:   R_Sharipov@ic.bashedu.ru
%           ruslan-sharipov@usa.net
% Web-page: http://www.geocities.com/CapeCanaveral/Lab/5341
%
\nopagenumbers
\catcode`@=11
\redefine\logo@{\relax}
\catcode`\@=\active
\def\startpage#1{\pageno=#1}
\pagewidth{12.8cm}
\pageheight{20cm}
\def\tr{\operatorname{tr}}
\def\const{\operatorname{const}}
\def\compos{\,\raise 1pt\hbox{$\sssize\circ$} \,}
\accentedsymbol\DotTildePhi{\Dot{\Tilde\varphi}}
\accentedsymbol\tx{\tilde x}
\Monograph
\topmatter
\title\chapter{5}
Normal shift in Riemannian manifolds.
\endtitle
\endtopmatter
\loadbold
\startpage{83}
\document
\head
\S\,1. Geodesic normal shift.
\endhead
     In this Chapter we return to the case of Riemannian manifolds
and derive main results of thesis. They consist in developing the
theory of {\bf Newtonian dynamical systems admitting the normal
shift} on such manifolds. As a starting point for this theory we
choose well-known classical construction of {\bf geodesic normal
shift}.\par
     Let $M$ be Riemannian manifold with metric tensor $\bold g$.
Tensor field $\bold g$ defines scalar product of tangent vectors
on the manifold $M$. Let's denote it as
$$
\bold g(\bold X,\bold Y)=(\bold X\,|\,\bold Y).\hskip -2em
\tag1.1
$$
Here $\bold X$ and $\bold Y$ are two vectors at one point $p\in M$.
Metric $\bold g$ defines metric connection $\Gamma$ with identically
zero torsion field $\bold T=0$. Its components in local map
$(U,x^1,\ldots,x^n)$ are determined by components of metric tensor
$\bold g$ according to the following standard formula (see \cite{12},
\cite{32}, \cite{76}, or \cite{77}):
$$
\Gamma^k_{ij}=\frac{1}{2}\sum^n_{s=1}g^{ks}\left(\frac{\partial g_{sj}}
{\partial x^i}+\frac{\partial g_{is}}{\partial x^j}-\frac{\partial g_{ij}}
{\partial x^s}\right),\hskip -2em
\tag1.2
$$
Quantities $g_{ij}$ and $\Gamma^k_{ij}$ in \thetag{1.2} do not depend
on the components of velocity vector  $v^1,\,\ldots,\,v^n$, i\.~e\. we
have ordinary (not extended) tensor field $\bold g$ and ordinary (not
extended) connection $\Gamma$. Applying velocity gradient $\tilde\nabla$
to the field $\bold g$, we get zero, while space gradient $\nabla$
applied to $\bold g$ coincides with traditional covariant differential of
this field. Therefore we have the equalities 
$$
\xalignat 2
&\nabla_\bold g=0,&&\tilde\nabla\bold g=0.\hskip -2em
\tag1.3
\endxalignat
$$
First equality in \thetag{1.2} is standard (see \cite{12}, \cite{32},
\cite{76}, or \cite{77}). It expresses the {\bf concordance of} metric
and connection).\par
    {\bf Geodesic lines} play important role in geometry. In local
coordinates  the are defined by ordinary differential equations
$$
\qquad\ddot x^k+\sum^n_{i=1}\sum^n_{j=1}\Gamma^k_{ij}\,\dot x^i\,
\dot x^j=0\text{, \ \ where\ \ \ }k=1,\,\ldots,\,n,
\hskip -2em
\tag1.4
$$
provided special parametrization is used, where $t$ is a length
of the curve measured from some fixed reference point on it. Let's
write the equations \thetag{1.4} in terms of covariant differentiation
with respect to parameter $t$ along the curve (such covariant
differentiations were considered in Chapter
\uppercase\expandafter{\romannumeral 4}):
$$
\xalignat 2
&\dot x^k=v^k,&&\nabla_tv^k=0.\hskip -2em
\tag1.5
\endxalignat
$$
The equations \thetag{1.5} form the system of $2n$ ordinary
differential equations of the first order. In setting up Cauchy
problem for these equations initial values of local coordinates
$x^1,\,\ldots,\,x^n$ and components of velocity vector
$v^1,\,\ldots,\,v^n$ are given:
$$
\xalignat 2
&\quad x^k(t)\,\hbox{\vrule height 8pt depth 8pt width 0.5pt}_{\,t=0}
=x^k(0),
&&v^k(t)\,\hbox{\vrule height 8pt depth 8pt width 0.5pt}_{\,t=0}=v^k(0).
\hskip -2em
\tag1.6
\endxalignat
$$
This means that we fix initial point $p=p(0)$ and initial vector
of velocity $\bold v(0)$ at the point $p(0)$. Note that equations
\thetag{1.4} are invariant under the following transformation of
independent variable: $t\to\const\cdot\,t$. Therefore without loss
of generality we can assume the initial vector $\bold v(0)$ in
\thetag{1.6} to be of unit length: $|\bold v(0)|=1$. Once fulfilled
for $t=0$ this condition remains fulfilled for nonzero $t$ as well:
$$
|\bold v(t)|=1.\hskip -2em
\tag1.7
$$
This is the very case when parameter $t$ is a length of curve
measured from some fixed reference point, it is called {\bf natural
parameter} on geodesic line.\par
\parshape 19 6.0cm 6.8cm 6.0cm 6.8cm 6.0cm 6.8cm 6.0cm 6.8cm
6.0cm 6.8cm 6.0cm 6.8cm 6.0cm 6.8cm 6.0cm 6.8cm 6.0cm 6.8cm
6.0cm 6.8cm 6.0cm 6.8cm 6.0cm 6.8cm 6.0cm 6.8cm 6.0cm 6.8cm
6.0cm 6.8cm 6.0cm 6.8cm 6.0cm 6.8cm 6.0cm 6.8cm 0cm 12.8cm
     Let $S$ be some hypersurface in $M$. For the sake of simplicity
$S$ assumed to be connected and simply connected smooth submanifold
of codimension $1$ without boundary such that the closure of $S$
is compact. \vadjust{\vskip -1.5cm\special{em:graph edis-05a.gif}
\vskip 1.5cm}Moreover, we assume that $S$ is compactly imbedded into
some other open hypersurface, i\.~e\. $S\Subset S'$. This strong
requirement eliminates most difficulties that could arise on the
boundary $\partial S$. Under the above assumptions hypersurface $S$
is orientable, and one can choose smooth field of unitary normal
vectors on it, which has smooth continuation to any point on the
boundary $\partial S$. This field defines unitary normal vector
$\bold n(p)$ at each point $p\in S$. Let's choose components of
$\bold n(p)$ as the components of initial velocity in setting up
initial data for Cauchy problem \thetag{1.6}. Then the condition
$|\bold v(0)|=1$ is fulfilled due to $|\bold n(p)|=1$. Let's write
\thetag{1.6} as
$$
\pagebreak
\xalignat 2
&\quad x^k(t)\,\hbox{\vrule height 8pt depth 8pt width 0.5pt}_{\,t=0}
=x^k(p),
&&v^k(t)\,\hbox{\vrule height 8pt depth 8pt width 0.5pt}_{\,t=0}=n^k(p).
\hskip -2em
\tag1.8
\endxalignat
$$
Solution of the equations \thetag{1.5} satisfying initial data \thetag{1.8}
defines the family of geodesic lines beginning at the points of hypersurface
$S$ and being perpendicular to $S$ at that points (See Fig\. 1.1). The value
of parameter $t$ on these lines coincides with their length measured from
starting points on $S$.\par
     Let $p_0$ be a pint on the hypersurface $S$, and let $U$ be some map
on $S$ covering this point. Suppose that $u^1,\,\ldots,\,u^{n-1}$ are local
coordinates of the point $p$ in the map $U$. Then geodesic lines obtained
as the solution of Cauchy problem \thetag{1.8} for the system of equations
\thetag{1.5} are represented by the set of $n$ functions
$$
\align
x^1&=x^1(u^1,\ldots,u^{n-1},t),\hskip -2em\\
.\ .\ &.\ .\ .\ .\ .\ .\ .\ .\ .\ .\ .\ .\ .\ .\ .\ .\ .
\hskip -2em\tag1.9\\
x^n&=x^n(u^1,\ldots,u^{n-1},t),\hskip -2em
\endalign
$$
where $x^1,\,\ldots,\,x^n$ are local coordinates in $M$. According
to well-known ``{existence and uniqueness}'' theorem (see \cite{82} or
\cite{83}) the domain of functions \thetag{1.9} is such that they
define the map $U'\times I_\delta\to M$, where $U'$ is an open
subset of local map $U$ containing the point $p_0$, and $I_\delta=
(-\delta,+\delta)$ is open $\delta$-neighborhood of zero on real
axis $\Bbb R$. Due to the compactness of the closure of $S$, or, being
more exact, due to $S\Subset S'$ local maps $U'\times I_\delta\to M$
can be glued into one global map $S\times I_\varepsilon\to M$.\par
     Let's differentiate the functions \thetag{1.9} with respect to
parameters $u^1,\,\ldots,\,u^{n-1}$, and let's form the following
vectors by the obtained derivatives:
$$
\bold K_i=\sum^n_{j=1}\frac{\partial x^j}{\partial u^i}\,
\frac{\partial}{\partial x^j}\text{, \ where \ }i=1,\,\ldots,\,n-1.
\hskip -3em
\tag1.10
$$
Differentiating the functions \thetag{1.9} with respect to their last
argument $t$, we form one more vector $\bold N$. It is defined by
the formula similar to \thetag{1.10}:
$$
\bold N=\sum^n_{j=1}\frac{\partial x^j}{\partial t}\,
\frac{\partial}{\partial x^j}.\hskip -3em
\tag1.11
$$
For $t=0$ vectors \thetag{1.10} coincide with coordinate tangent
vectors to $S$, while vector \thetag{1.11} coincides with normal
vector of $S$. Therefore vectors $\bold K_1,\,\ldots,\,\bold K_{n-1},
\,\bold N$ are linearly independent for $t=0$. If we denote $t=u^n$,
and if we build Jacoby matrix
$$
J=\Vmatrix
\dsize\frac{\partial x^1}{\partial u^1} &\hdots
&\dsize\frac{\partial x^1}{\partial u^n}\\
\vspace{1ex}
\vdots & \ddots & \vdots\\
\vspace{2ex}
\dsize\frac{\partial x^n}{\partial u^1} &\hdots
&\dsize\frac{\partial x^n}{\partial u^n}
\endVmatrix,
$$
it will be non-degenerate at the point $p_0$. This means, that the
variables $u^1,\,\ldots,\,u^n$ (where $u^n=t$) could be used as local
coordinates on $M$ in some neighborhood of the point $p_0$. These
coordinates are called {\bf semigeodesic coordinates} (see \cite{**}).
Hypersurface $S$ in these coordinates is described by the equation
$u^n=0$.\par
    Let's consider the equation $u^n=t=\const$ in semigeodesic coordinates.
It determines a piece of some hypersurface. Due to $S\Subset S'$ we can
find $\varepsilon>0$ such that for $t\in I_\varepsilon$ various pieces
of hypersurfaces $u^n=t$ are glued into an open hypersurface $S_t$
diffeomorphic to $S$. At the expense of choosing sufficiently small value
of $\varepsilon$ we can reach the situation, when for $t\in I_\varepsilon$
none of hypersurfaces $S_t$ will intersect another, nor will it have
self-intersection. Therefore $S$ has some neighborhood in $M$ diffeomorphic
to $S\times I_\varepsilon$, and for each $t\in I_\varepsilon$ there is
a diffeomorphism $f_t\!:\,S\to S_t$. Such diffeomorphism built by geodesic
lines as described above is called {\bf geodesic shift} of hypersurface
$S$ to the distance $t$. It's important  to  emphasize  that  any 
hypersurface
$S_t$ has compact closure, and there exists some open hypersurface $S'_t$
such that $S_t\Subset S'_t$.\par
    Let's study the geodesic shift of hypersurface $S$ in semigeodesic
coordinates $x^1=u^1,\,\ldots,\,x^n=u^n$ associated with it. In this case
$x^n=t$ is a parameter of shift, and vectors \thetag{1.10} are coordinate
tangent vectors to hypersurfaces $S_t$ described by the equations $x^n=t
=\const$. Due to $x^j=u^j$ they are the followings:
$$
\bold K_1=\frac{\partial}{\partial x^1},\ \ldots,\ \bold K_{n-1}
=\frac{\partial}{\partial x^{n-1}}.\hskip -3em
\tag1.12
$$
Geodesic lines implementing geodesic shift of hypersurface $x^n=0$
in semigeodesic coordinates are described by functions
$$
x^1(t)=\const,\ \ldots,\ x^{n-1}(t)=\const,\ x^n(t)=t.\hskip -3em
\tag1.13
$$
Let's substitute \thetag{1.13} into the equations of geodesic lines
\thetag{1.4}. This yields
$$
\Gamma^k_{nn}=0\text{\ \ for all \ }k=1,\,\ldots,\,n.\hskip -3em
\tag1.14
$$
Having calculated connection components $\Gamma^k_{nn}$ by formula
\thetag{1.2}, due to non-degene\-racy of the matrix of metric tensor
$g^{ks}$ from \thetag{1.14} we derive the relationships 
$$
2\,\frac{\partial g_{sn}}{\partial x^n}-\frac{\partial g_{nn}}
{\partial x^s}=0\text{\ \ for all \ }s=1,\,\ldots,\,n.\hskip -3em
\tag1.15
$$
Particularly for $s=n$ from \thetag{1.15} we extract the following
equation:
$$
\frac{\partial g_{nn}}{\partial x^n}=0.\hskip -3em
\tag1.16
$$
Thus, diagonal component $g_{nn}$ in metric tensor doesn't depend on
shift parameter $x^n=t$; and $g_{nn}$ is the square of the length of
$n$-th coordinate vector, \pagebreak which in semigeodesic coordinates
coincides with vector \thetag{1.11}. The latter one plays the role of
velocity vector on geodesic lines used to shift $S$:
$$
\bold K_n=\bold N=\bold v(t)=\frac{\partial}{\partial x^n}.
\hskip -3em
\tag1.17
$$
For $x^n=t=0$ velocity vector $\bold v(t)$ coincides with unitary
normal vector to $S$ (see \thetag{1.8}). Hence $g_{nn}=1$ for
$x^n=0$. From this equality combined with \thetag{1.16} we derive
$$
g_{nn}\equiv 1.\hskip -3em
\tag1.18
$$
In particular, \thetag{1.18} could be the proof for the equality
\thetag{1.7} for the velocity vector on geodesic lines. Substituting
\thetag{1.18} into \thetag{1.15} for $s\neq n$, we
immediately get
$$
\frac{\partial g_{sn}}{\partial x^n}=0\text{\ \ for \ }s=1,\,\ldots,
\,n-1.\hskip -3em
\tag1.19
$$
Components $g_{sn}$ in metric tensor are the scalar products of vectors
\thetag{1.12} and \thetag{1.17}:
$$
g_{sn}=(\bold K_s\,|\,\bold K_n).\hskip -3em
\tag1.20
$$
At the initial instant, when $x^n=t=0$, scalar products \thetag{1.20}
are equal to zero, since $\bold K_n$ coincides with unitary normal
vector on $S$, while vectors \thetag{1.12} are tangent to $S$. From
this fact combined with the equations \thetag{1.19} we get
$$
g_{sn}\equiv 0\text{\ \ for all \ }s=1,\,\ldots,
\,n-1.\hskip -3em
\tag1.21
$$
The equalities \thetag{1.21} reflect very important geometric property
of geodesic shift. They mean that not only the initial
hypersurface $S$ is perpendicular to geodesic lines used to shift it,
but all other hypersurfaces $S_t$ generated by shift procedure are
perpendicular as well. Therefore geodesic shift is a {\bf normal shift}
of hypersurfaces.
\head
\S\,2. Normal shift along trajectories of dynamical system.
\endhead
     Geodesic normal shift described in \S\,1 is a well-known classical
construction. The {\bf main idea} that brought about the advent of the
theory of {\bf dynamical systems admitting the normal shift} consists in
replacing geodesic lines by {\bf wider class of parametric curves}.
As a first pretender to this role we considered the trajectories of
Newtonian dynamical systems, since
\roster
\rosteritemwd=5pt
\item they possess prescribed parametrization with the parameter having
      physical meaning of time;
\item they are associated with positively definite quadratic form of
      kinetic energy that equips configuration space with the structure
      of Riemannian manifold;
\item in local coordinates they are described by systems of differential
      equations compatible with Cauchy problems of the form \thetag{1.6}.
\endroster
     Let's consider Newtonian dynamical system with force field $\bold F$,
configuration space of which is a Riemannian manifold $M$ with metric
$\bold g$. Let's choose some local map $(U,x^1,\ldots,x^n)$ on $M$ and
canonically associated local coordinates on tangent bundle $TM$
(see \S\,6 in Chapter \uppercase\expandafter{\romannumeral 2}). Then
trajectories of Newtonian dynamical system with the field $\bold F$
can be constructed by solving the system of differential equations
\thetag{4.1} from Chapter \uppercase\expandafter{\romannumeral 2}.
By analogy with \thetag{1.5} these equations could be written in terms
of covariant derivative with respect to parameter $t$:
$$
\xalignat 2
&\dot x^k=v^k,&&\nabla_tv^k=F^k.\hskip -2em
\tag2.1
\endxalignat
$$
\definition{Definition 2.1} We say that hypersurface $S$ in $M$
belongs to {\bf transformation class}, if it is connected and simply
connected open hypersurface with compact closure such that
$S\Subset S'$ for some other open hypersurface $S'$.
\enddefinition
    Note that closed hypersurface without boundary belongs to
transformation class, if it is connected, simply connected, and
compact, since for closed hypersurface $\overline{S}=S$. Here
we can take $S'=S$.\par
    Suppose that $S\subset M$ is a hypersurface from transformation class.
Let's fix the field of unitary normal vectors $\bold n(p)$ on $S$, and
let's set up the following initial data for the system of differential
equations \thetag{2.1} at the points of $S$:
$$
\xalignat 2
&\quad x^k(t)\,\hbox{\vrule height 8pt depth 8pt width 0.5pt}_{\,t=0}
=x^k(p),
&&v^k(t)\,\hbox{\vrule height 8pt depth 8pt width 0.5pt}_{\,t=0}=
\nu(p)\cdot n^k(p).\hskip -2em
\tag2.2
\endxalignat
$$
Solution of Cauchy problem with initial data \thetag{2.2} determines
the family of trajectories of Newtonian dynamical system \thetag{2.1}
beginning at the points of $S$ and outgoing from $S$ along normal
vector of this hypersurface (see Fig\. 1.1). The value of $\nu(p)$
determines modulus of velocity vector at the initial point $p$ on such
trajectory. In the case of geodesic normal shift the value of $\nu(p)$
is chosen to be equal to unity everywhere on $S$. Here we shall not
keep this restriction, we shall require only that $\nu(p)\neq 0$ for
all points $p\in S$.
\definition{Definition 2.2} Let $S\subset M$ be a hypersurface from
transformation class. We say that smooth function $\nu(p)$ on $S$
belongs to {\bf transformation class}, if it has smooth expansion
to some bigger open hypersurface $S'$ such that $S\Subset S'$, and if
it is nonzero everywhere on the closure $\overline{S}$.
\enddefinition
    Suppose that on the hypersurface $S$ belonging to the transformation
class on $M$ we defined some function $\nu(p)$ belonging to the
transformation class on $S$. For each point $p=p_0$ on $S$ we draw the
trajectory of dynamical system \thetag{2.1} passing through $p_0$ and
satisfying initial data \thetag{2.2} at this point. Then consider the
point $p(t)=p(p_0,t)$ on such trajectory corresponding to the value $t$
of time variable, and associate it with $p_0$. So we have constructed
a transformation of the hypersurface $S$ that maps $p_0$ to $p(t)$. It
is called the {\bf shift along trajectories of dynamical system} for
the time $t$.\par
     In the case of shift along trajectories of Newtonian dynamical
system, as well as in the case of geodesic shift, we can point out
some interval $I_\varepsilon$ for time variable $t$ such that the
map $S\times I_\varepsilon\to M$ is injective and has no singular
points. This map stratifies some full-dimensional neighborhood of $S$
into a family of hypersurfaces $S_t$ diffeomorphic to $S$: for any
$t\in I_\varepsilon$ we have diffeomorphism $f_t\!:\,S\to S_t$.
For $t\in I_\varepsilon$ hypersurfaces $S_t$ all are in transformation
class, they have no self-intersection points and they do not intersect
each other.
\definition{Definition  2.3}  Shift  of  hypersurface  $S$  along 
trajectories
of Newtonian dynamical system \thetag{2.1} initiated by the function
$\nu(p)$ in Cauchy problem \thetag{2.2} is called {\bf normal shift}, if
for parameter $t$ in some interval $I_\varepsilon=(-\varepsilon,
+\varepsilon)$ all hypersurfaces $S_t$ are orthogonal to the trajectories
used to shift $S$.
\enddefinition
\noindent The definition~2.3 is quite necessary, since, in contrast to
geodesic normal shift, here the orthogonality of $S_t$ and trajectories
is not an unconditional property of shift. Presence (or absence) of this
property essentially depends on how we deal with arbitrariness in the
choice of the function $\nu(p)$ on $S$ and in the choice of force field
$\bold F$ for dynamical system.
\head
\S\, 3. Differential equations for the vector of variation.
\endhead
     Let's choose some local map $(U,x^1,\ldots,x^n)$ on $M$ and consider
the shift of some hypersurface $S$ along trajectories of Newtonian
dynamical system \thetag{2.1}. If $u^1,\,\ldots,\,u^{n-1}$ are local
coordinates on $S$, then trajectories of dynamical system \thetag{2.1}
beginning at the points of this hypersurface are expressed by functions
$$
\align
x^1&=x^1(u^1,\ldots,u^{n-1},t),\hskip -2em\\
.\ .\ &.\ .\ .\ .\ .\ .\ .\ .\ .\ .\ .\ .\ .\ .\ .\ .\ .
\hskip -2em\tag3.1\\
x^n&=x^n(u^1,\ldots,u^{n-1},t).\hskip -2em
\endalign
$$
Parameter $t$ in \thetag{3.1} is the main parameter, it has physical
meaning of time, other parameters $u^1,\,\ldots,\,u^{n-1}$ are auxiliary
ones. If we choose one of auxiliary parameters $u^i$ and fix all others,
then we get the situation considered in \S\,5 of Chapter
\uppercase\expandafter{\romannumeral 4}. In this situation partial
derivatives
$$
\frac{\partial x^1}{\partial u^i},\ \ldots,\ \frac{\partial x^n}
{\partial u^i}\hskip -2em
\tag3.2
$$
are the components of the vector of variation. We denote it by $\boldsymbol
\tau_i$. On the other hand, these derivatives are the components of the
vector $\bold K_i$ from \thetag{1.10} as well. Thus, on the trajectories
of dynamical system \thetag{2.1} used in the construction of shift we have
several vectors of variation simultaneously:
$$
\boldsymbol\tau_1=\bold K_1,\ \ldots,\ \boldsymbol\tau_{n-1}=\bold
K_{n-1}.\hskip -2em
\tag3.3
$$
From geometrical point of view vectors \thetag{3.3} are
interpreted as coordinate tangent vectors on hypersurface $S_t$ in
local coordinates $u^1,\,\ldots,\,u^{n-1}$ induced to $S_t$ from $S$
by means of shift diffeomorphism $f_t\!:\,S\to S_t$.\par
     Let's study time dependence of the vectors \thetag{3.3}, i\.~e\.
their evolution as we move along trajectories of the system
\thetag{2.1}. Components of any of them appear to satisfy the same
differential equations, which can be obtained by linearization of
the equations \thetag{2.1}. In deriving these equations we denote
$u^i=u$ and $\bold K_i=\boldsymbol\tau_i=\boldsymbol\tau$ for to
simplify further calculations. The equations \thetag{2.1} determine
not only the functions $x^i(u^1,\ldots,u^{n-1},t)$ in \thetag{3.1},
but the functions
$$
\align
v^1&=v^1(u^1,\ldots,u^{n-1},t),\hskip -2em\\
.\ .\ &.\ .\ .\ .\ .\ .\ .\ .\ .\ .\ .\ .\ .\ .\ .\ .\ .
\hskip -2em\tag3.4\\
v^n&=v^n(u^1,\ldots,u^{n-1},t)\hskip -2em
\endalign
$$
as well. Therefore we have natural lift of trajectories to the
tangent bundle $TM$ (see definition~2.1 in Chapter
\uppercase\expandafter{\romannumeral 4}). Let's apply the operator
$\nabla_u$ of covariant differentiation with respect to parameter
$u$ (see \S\,5 in Chapter \uppercase\expandafter{\romannumeral 4})
to both sides of second part of the equations \thetag{2.1}. The
result can be written as a vectorial equality
$$
\nabla_u\nabla_t\bold v=\nabla_u(\bold F).\hskip -2em
\tag3.5
$$
In order to transform left hand side of \thetag{3.5} we transpose
covariant differentiations $\nabla_u$ and $\nabla_t$ by means of
formula \thetag{5.14} from Chapter \uppercase\expandafter{\romannumeral
4}:
$$
\nabla_u\nabla_t\bold v=[\nabla_u,\nabla_t]\bold v+\nabla_t\nabla_u
\bold v=\bold S\bold v+\nabla_t\nabla_u\bold v.\hskip -2em
\tag3.6
$$
Operator field $\bold S$ in \thetag{3.6} is determined by formula
\thetag{5.1} from Chapter \uppercase\expandafter{\romannumeral 4}.
Applying this formula, we take into account that dynamic curvature
field $\bold D$ and torsion field $T$ for metric connection
\thetag{1.2} are zero. Therefore
$$
\nabla_u\nabla_t\bold v=\nabla_t\nabla_u\bold v+\bold R(\boldsymbol
\tau,\bold v)\bold v.\hskip -2em
\tag3.7
$$
In order to calculate $\nabla_u\bold v$ we use formula \thetag{5.13}
from Chapter \uppercase\expandafter{\romannumeral 4} and remember that
torsion is zero. Then from \thetag{3.7} we derive
$$
\nabla_u\nabla_t\bold v=\nabla_{tt}\boldsymbol\tau+\bold R(\boldsymbol
\tau,\bold v)\bold v.\hskip -2em
\tag3.8
$$
In order to transform right hand side of equality \thetag{3.5} we use
formula \thetag{5.19} or formula \thetag{5.20} from Chapter
\uppercase\expandafter{\romannumeral 4} and remember that $\bold T=0$:
$$
\nabla_u(\bold F)=C(\boldsymbol\tau\otimes\nabla\bold F)+
C(\nabla_t\boldsymbol\tau\otimes\tilde\nabla\bold F).
\tag3.9
$$
Now, equating the expressions \thetag{3.8} and \thetag{3.9}, we
obtain vectorial differential equation for the vector of variation
$\boldsymbol\tau$:
$$
\pagebreak
\nabla_{tt}\boldsymbol\tau=C(\nabla_t\boldsymbol\tau\otimes
\tilde\nabla\bold F)+C(\boldsymbol\tau\otimes\nabla\bold F)
-\bold R(\boldsymbol\tau,\bold v)\bold v.\hskip -3em
\tag3.10
$$
Written in terms of components in some local map, \thetag{3.10}
looks like 
$$
\nabla_{tt}\tau^k=-\sum^n_{m=1}\sum^n_{i=1}\sum^n_{j=1}R^k_{mij}\,
\tau^i\,v^j\,v^m+\sum^n_{m=1}\nabla_t\tau^m\,\tilde\nabla_mF^k
+\sum^n_{m=1}\tau^m\,\nabla_mF^k.
\hskip -2em
\tag3.11
$$
From \thetag{3.11} we conclude that components of variation vector
$\boldsymbol\tau$ satisfy the system of $n$ linear ordinary equation
of second order respective to time variable $t$. Coefficients in these
linear equations are depending on $t$. They are determined by the force
field $\bold F$ of dynamical system, and they depend on the choice of its
particular trajectory, i\.~e\. they depend on the functions \thetag{3.1}
and \thetag{3.4}.
\head
\S\,4. Function of deviation and its derivatives.
\endhead
     As a vector of variation $\boldsymbol\tau$ in the equations
\thetag{3.10} we can choose any one of vectors $\boldsymbol\tau_1,\,
\ldots,\,\boldsymbol\tau_{n-1}$ from \thetag{3.3}. In the case, when
$f_t\!:\,S\to S_t$ is a normal shift, they all are perpendicular to
the vector $\bold v(t)$, which is tangent to the trajectory. Therefore
as a measure of deviation from normality we use the following scalar
products:
$$
\varphi_i=(\boldsymbol\tau_i\,|\,\bold v)\text{, \ \ where \ \ \ }
i=1,\,\ldots,\,n-1.\hskip -2em
\tag4.1
$$
In the situation of normal shift all these functions are identically
zero for all trajectories of shift.
\definition{Definition 4.1} Scalar product of velocity vector
$\bold v(t)$ with the vector of variation $\boldsymbol\tau(t)$
is called the function of {\bf deviation} on the trajectory of
Newtonian dynamical system.
\enddefinition
    Suppose that $\varphi(t)$ is a function of deviation on some
particular trajectory of Newtonian dynamical system:
$$
\varphi=(\bold v\,|\,\boldsymbol\tau).\hskip -2em
\tag4.2
$$
Let's calculate two time derivatives of this function $\dot\varphi(t)$
and $\ddot\varphi(t)$. For the first derivative of the function
\thetag{4.2} we have
$$
\dot\varphi=\nabla_t\varphi=\nabla_t(\bold v\,|\,\boldsymbol\tau)=
(\nabla_t\bold v\,|\,\boldsymbol\tau)+(\bold v\,|\,\nabla_t\boldsymbol
\tau)\hskip -2em
\tag4.3
$$
Differentiating scalar product $(\bold v\,|\,\boldsymbol\tau)$ in
\thetag{4.3}, we used Leibniz rule. This is correct, since we
can do the following more detailed calculations
$$
\gathered
\nabla_t(\bold v\,|\,\boldsymbol\tau)=\nabla_tC(\bold g\otimes\bold v
\otimes\boldsymbol\tau)=C(\nabla_t\bold g\otimes\bold v\otimes\boldsymbol
\tau)+C(\bold g\otimes\nabla_t\bold v\otimes\boldsymbol\tau)\,+\\
\vspace{1ex}
+\,C(\bold g\otimes\bold v\otimes\nabla_t\boldsymbol\tau)=(\nabla_t\bold
v\,|\,\boldsymbol\tau)+(\bold v\,|\,\nabla_t\boldsymbol\tau)
+C(\nabla_t\bold g\otimes\bold v\otimes\boldsymbol\tau).
\endgathered\hskip -2em
\tag4.4
$$
Covariant derivative $\nabla_t\bold g$ in \thetag{4.4} is equal to zero.
\pagebreak
This follows from the formula \thetag{4.3} in Chapter
\uppercase\expandafter{\romannumeral 4} and the relationships \thetag{1.3}
that express concordance of metric and connection. Indeed, we have
$$
\nabla_t\bold g=C(\bold v\otimes\nabla\bold g)+C(\nabla_t\bold v\otimes
\tilde\nabla\bold g)=0.\hskip -2em
\tag4.5
$$
For the further transformation of \thetag{4.3} we use the equations
of dynamics \thetag{2.1}, second part of which can be written in vectorial
form:
$$
\nabla_t\bold v=\bold F.\hskip -2em
\tag4.6
$$
Substituting \thetag{4.6} into the formula \thetag{4.3}, for the
derivative $\dot\varphi$ we get
$$
\dot\varphi=(\bold F\,|\,\boldsymbol\tau)+(\bold v\,|\,
\nabla_t\boldsymbol\tau).\hskip -2em
\tag4.7
$$\par
    Now let's calculate second derivative $\ddot\varphi$. In order to
do it let's differentiate \thetag{4.7}, taking into account the
equation \thetag{4.6} thereby:
$$
\ddot\varphi=\nabla_t\dot\varphi=(\nabla_t\bold F\,|\,\boldsymbol\tau)
+2\,(\bold F\,|\,\nabla_t\boldsymbol\tau)+(\bold v\,|\,\nabla_{tt}
\boldsymbol\tau).\hskip -2em
\tag4.8
$$
For to calculate $\nabla_t\bold F$ we apply formula \thetag{4.3}
from Chapter \uppercase\expandafter{\romannumeral 4}:
$$
\nabla_t\bold F=C(\bold v\otimes\nabla\bold F)+
C(\bold F\otimes\tilde\nabla\bold F).\hskip -2em
\tag4.9
$$
Second derivative of variation vector $\nabla_{tt}\boldsymbol\tau$
is expressed through $\boldsymbol\tau$ and through first derivative
$\nabla_t\boldsymbol\tau$ by means of the equation \thetag{3.10}.
If additionally we take into account \thetag{4.9}, then for the
derivative $\ddot\varphi$ we get
$$
\gathered
\ddot\varphi=(C(\bold v\otimes\nabla\bold F)\,|\,\boldsymbol\tau)+
(C(\bold F\otimes\tilde\nabla\bold F)\,|\,\boldsymbol\tau)
+2\,(\bold F\,|\,\nabla_t\boldsymbol\tau)\,+\\
\vspace{1ex}
+\,(\bold v\,|\,C(\boldsymbol\tau\otimes\nabla\bold F))
+(\bold v\,|\,C(\nabla_t\boldsymbol\tau\otimes\tilde\nabla\bold F))
-(\bold v\,|\,\bold R(\boldsymbol\tau,\bold v)\bold v).
\endgathered\hskip -2em
\tag4.10
$$
Last term in \thetag{4.10} is equal to zero. This follows from some
properties of curvature tensor of metric connection \thetag{1.2}.
Indeed, let's write the expression $(\bold v\,|\,\bold R(\boldsymbol
\tau,\bold v)\bold v)$ in local coordinates. Here we have
$$
(\bold v\,|\,\bold R(\boldsymbol\tau,\bold v)\bold v)=
\sum^n_{k=1}\sum^n_{m=1}\sum^n_{i=1}\sum^n_{j=1}R_{kmij}\,
v^k\,\tau^i\,v^j\,v^m=0.\hskip -3em
\tag4.11
$$
Sum in \thetag{4.11} vanishes due to skew symmetry of $R_{kmij}$
with respect to first two indices $k$ and $m$ (see \cite{12}, \cite{32},
\cite{76}, or \cite{77}). Therefore
$$
\gathered
\ddot\varphi=(C(\bold v\otimes\nabla\bold F)\,|\,\boldsymbol\tau)+
(C(\bold F\otimes\tilde\nabla\bold F)\,|\,\boldsymbol\tau)\,+\\
\vspace{1ex}
+\,2\,(\bold F\,|\,\nabla_t\boldsymbol\tau)
+(\bold v\,|\,C(\boldsymbol\tau\otimes\nabla\bold F))+
(\bold v\,|\,C(\nabla_t\boldsymbol\tau\otimes\tilde\nabla\bold F)).
\endgathered\hskip -2em
\tag4.12
$$\par
    Let's write the equalities \thetag{4.2}, \thetag{4.7}, and
\thetag{4.12} in local coordinates. For the function of deviation
$\varphi$ and for its first derivative $\dot\varphi$ we have
$$
\align
\varphi&=\sum^n_{i=1}v_i\,\tau^i,\hskip -3em
\tag4.13\\
\dot\varphi&=\sum^n_{i=1}F_i\,\tau^i+\sum^n_{i=1}v_i\,\nabla_t\tau^i
\hskip -3em
\tag4.14
\endalign
$$
For the sake of brevity in \thetag{4.13} and \thetag{4.14} we used
covariant components of vectors  $\bold v$ and $\bold F$, which are
obtained by lowering the indices:
$$
\xalignat 2
&v_i=\sum^n_{j=1}g_{ij}\,v^j,
&&F_i=\sum^n_{j=1}g_{ij}\,F^j.
\hskip -3em
\tag4.15
\endxalignat
$$
Procedure of lowering indices \thetag{4.15} is commutating with
various covariant differentiations due to the concordance of metric
and connection expressed by \thetag{1.3} (see more details in
\cite{77}). This procedure is invertible, therefore we can make
almost no difference between covariant and contravariant components
of tensorial objects in Riemannian geometry.\par
    Now let's write in local coordinates the expression \thetag{4.12}
for the second derivative of the function of variation:
$$
\gathered
\ddot\varphi=\sum^n_{i=1}\left(2\,F_i+\shave{\sum^n_{j=1}} v^j\,
\tilde\nabla_iF_j\right)\nabla_t\tau^i\,+\hskip -2em\\
+\,\sum^n_{i=1}\left(\,\shave{\sum^n_{j=1}}v^j\left(\nabla_jF_i+
\nabla_iF_j\right)+\shave{\sum^n_{j=1}} F^j\,\tilde\nabla_jF_i
\right)\tau^i.\hskip -2em
\endgathered
\tag4.16
$$
In this formula \thetag{4.16} we use covariant and contravariant
components of force field $\bold F$ simultaneously.
\head
\S\,5. Projectors defined by the vector of velocity.
\endhead
    Function of deviation $\varphi$ depends only on orthogonal
projection of variation vector $\boldsymbol\tau$ to the direction
of velocity vector. This follows from \thetag{4.13} or from the
initial formula \thetag{4.2} for the function $\varphi$. Similarly,
first derivative $\dot\varphi$ also depends only on projection of
variation vector $\boldsymbol\tau$ to the direction of velocity
vector. In what follows it will be important to subdivide all entries
of $\boldsymbol\tau$ and $\nabla_t\boldsymbol\tau$ in formulas for
$\varphi$, $\dot\varphi$, and $\ddot\varphi$ into two parts, one
directed along the velocity vector $\bold v$, and the second directed
perpendicular to $\bold v$. With this aim in mind we define a few
auxiliary tensor fields from extended algebra $\bold T(M)$.\par
     Scalar field of modulus of velocity vector $v=|\bold v|$ is an
extended tensor field of type $(0,0)$. \pagebreak Its domain is the
whole tangent bundle $TM$, but it is smooth only at that points
$q=(p,\bold v)$, where $\bold v\neq 0$. In local coordinates for the
field $v$ we have
$$
v=\sqrt{\shave{\sum^n_{i=1}}\shave{\sum^n_{j=1}}g_{ij}\,v^i\,v^j}.
\hskip -2em
\tag5.1
$$
Denote by $\bold N$ the field of unit vectors directed along the velocity
vector $\bold v$. Then
$$
\bold N=\frac{\bold v}{v}.\hskip -2em
\tag5.2
$$
From \thetag{5.2} we see that $\bold N$ is an extended vector field,
the domain of which is the whole tangent bundle $TM$ except for those
points $q=(p,\bold v)$, where $\bold v=0$. The same condition $\bold v
\neq 0$ restricts domain of the following two extended operator fields:
$\bold Q$ and $\bold P$. $\bold Q$ is an operator of orthogonal projection
to the direction of the vector $\bold v$. Operator $\bold P$ is a projector
onto the hyperplane perpendicular to $\bold v$. Then
$$
\bold P+\bold Q=\bold 1.
\tag5.3
$$
Here $\bold 1$ is the field of identical operators. In local coordinates
operator fields $\bold P$ and $\bold Q$ have the following components:
$$
\xalignat 2
&\quad P^i_j=\delta^i_j-N^i\,N_j,
&&Q^i_j=N^i\,N_j.\hskip -2em
\tag5.4
\endxalignat
$$
Here $N^i$ and $N_j$ are the covariant and contravariant components
of the field \thetag{5.2}.\par
    From \thetag{5.1}, \thetag{5.2}, and \thetag{5.4} one can easily
derive the formulas for differentiating all the above defined
auxiliary fields. Let's write these formulas in local coordinates.
In case of vector field of velocity $\bold v$ and scalar field
$v=|\bold v|$ we have
$$
\xalignat 2
\nabla_kv^i&=0,
&\tilde\nabla_kv^i&=\delta^i_k,\hskip -2em
\tag5.5\\
\nabla_kv&=0,
&\tilde\nabla_kv&=N_k.\hskip -2em
\tag5.6
\endxalignat
$$
In case of vector field $\bold N$  from the relationships \thetag{5.5} and
\thetag{5.6} we derive
$$
\xalignat 2
\nabla_kN^i&=0,&\tilde\nabla_kN^i&=\frac{1}{v}\,P^k_i.\hskip -2em
\tag5.7
\endxalignat
$$
Spatial gradients of projector fields $\bold P$ and $\bold Q$ are
identically zero:
$$
\xalignat 2
\nabla_kP^i_j&=0,&\nabla_kQ^i_j&=0.\hskip -2em
\tag5.8
\endxalignat
$$
For velocity gradients of these two fields we have the formulas:
$$
\aligned
\tilde\nabla_kP^i_j=-&\frac{N_j\,P^i_k}{v}
-\sum^n_{r=1}\frac{g_{jr}\,P^r_k\,N^i}{v},\\
\vspace{0.5ex}
\tilde\nabla_kQ^i_j=&\frac{N_j\,P^i_k}{v}+\sum^n_{r=1}
\frac{g_{jr}\,P^r_k\,N^i}{v}.
\endaligned\hskip -2em
\tag5.9
$$
Spatial gradients are identically zero not only for the fields
$\bold P$ and $\bold Q$; for the fields $v$ and $\bold N$ they vanish as
well. This reflects the fact that they all are produced by the velocity
field $\bold v$ and Riemannian metric field $\bold g$, which satisfy the
relationships \thetag{1.3}. Derivation of formulas \thetag{5.5},
\thetag{5.6}, \thetag{5.7}, \thetag{5.8}, and \thetag{5.9} is rather
simple. It is based on formulas \thetag{7.3} and \thetag{7.4} from 
Chapter \uppercase\expandafter{\romannumeral 3}. We shall not give it
here.\par
    Let's use the relationship \thetag{5.3} in order to transform
first summand in formula \thetag{4.14} for the derivative $\dot\varphi$.
In local coordinates \thetag{5.3} is written as $\delta^i_j=P^i_j+Q^i_j$.
We can insert $\delta^i_j$ into \thetag{4.14} and add one more summation
there:
$$
\dot\varphi=\sum^n_{i=1}\sum^n_{j=1}F_i\,\delta^i_j\tau^j
+\sum^n_{i=1}v_i\,\nabla_t\tau^i\hskip -3em
\tag5.10
$$
Now, substituting $P^i_j+Q^i_j$ for $\delta^i_j$ in \thetag{5.10},
we get the formula
$$
\dot\varphi=\sum^n_{i=1}\sum^n_{j=1}F_i\,P^i_j\tau^j+
\sum^n_{i=1}\sum^n_{j=1}\frac{F_i\,N^i}{v}\,v_j\tau^j
+\sum^n_{i=1}v_i\,\nabla_t\tau^i.\hskip -3em
\tag5.11
$$
Formula \thetag{5.11} can be written in coordinate-free form as follows:
$$
\dot\varphi=(\bold F\,|\,\bold P\,\boldsymbol\tau)+
\frac{(\bold F\,|\,\bold N)}{v}\,(\bold v\,|\,\boldsymbol\tau)+
(\bold v\,|\,\nabla_t\boldsymbol\tau).\hskip -3em
\tag5.12
$$
Last two terms in formula \thetag{5.12} contain scalar products of
$\boldsymbol\tau$ and $\nabla_t\boldsymbol\tau$ with the vector of
velocity, while the first term does not depend on the projection of
these two vectors to the direction of $\bold v$.\par
    Now let's perform some similar transformations with the formula
\thetag{4.16} for the second derivative $\ddot\varphi$. Here it's
convenient to introduce auxiliary covectorial fields $\boldsymbol\alpha$
and $\boldsymbol\beta$ with the following components
$$
\aligned
\alpha_i&=2\,F_i+\sum^n_{j=1} v^j\,\tilde\nabla_iF_j,\hskip -2em\\
\beta_i&=\sum^n_{j=1}v^j\left(\nabla_jF_i+
\nabla_iF_j\right)+\sum^n_{j=1} F^j\,\tilde\nabla_jF_i.\hskip -2em
\endaligned
\tag5.13
$$
Then for the second derivative of the function of deviation we get
$$
\aligned
\ddot\varphi&=\boldsymbol\alpha(\bold P\nabla_t\boldsymbol\tau)
+\boldsymbol\beta(\bold P\boldsymbol\tau)
+\frac{\boldsymbol\alpha(\bold N)}{v}\,(\bold v\,|\,\nabla_t\boldsymbol\tau)
+\frac{\boldsymbol\beta(\bold N)}{v}\,(\bold v\,|\,\boldsymbol\tau).
\endaligned\hskip -2em
\tag5.14
$$
Last two terms in \thetag{5.14} contain scalar product of $\boldsymbol\tau$
and $\nabla_t\boldsymbol\tau$ with the vector of velocity $\bold v$. First
two terms do not depend oh the projections of $\boldsymbol\tau$ and
$\nabla_t\boldsymbol\tau$ to the direction of the vector $\bold v$.
\head
\S\,6. The condition of weak normality.
\endhead
     Let's look on the equations \thetag{3.11} and on the formulas
\thetag{4.13}, \thetag{4.14}, and \thetag{4.16} from the point of view
of differential equations theorist. When we fix some trajectory of
dynamical system \thetag{2.1}, the equations \thetag{3.11} form the
system of linear ordinary differential equations, overall order of
which is equal to $2n$. Denote by $\goth T$ the space of solutions of
this system of equations. Then $\dim\goth T=2n$. As the coordinates in
this space we can use components of two vectors $\boldsymbol\tau(t)$
and $\nabla_t\boldsymbol\tau(t)$ taken for some fixed value of
parameter $t=t_0$. Comparing formulas \thetag{4.13}, \thetag{4.14},
and \thetag{4.16}, we see that the function of deviation $\varphi$
and its derivatives $\dot\varphi$ and $\ddot\varphi$ depend linearly
on $\boldsymbol\tau$ and $\nabla_t\boldsymbol\tau$. This is true
for all time derivatives of the function of deviation, it is proved
by applying repeatedly $\nabla_t$ to the equality \thetag{4.12}.
So the value of $\varphi$ by itself, and values of all its
derivatives (corresponding to some fixed value of $t=t_0$)
$$
\varphi(t_0),\ \dot\varphi(t_0),\ \ddot\varphi(t_0),\
\varphi^{(3)}(t_0),\ \ldots,\ \varphi^{(2n)}(t_0),
\hskip -2em
\tag6.1
$$
can be considered as {\bf linear functionals} on the space
$\goth T$. The number of functionals \thetag{6.1} is equal to
$2n+1$, it is one as more than the dimension of dual space
$\goth T^*$. Therefore these functionals are linearly dependent.
Hence 
$$
\sum^{2n}_{i=0}C_i(t_0)\,\varphi^{(i)}(t_0)=0.\hskip -2em
\tag6.2
$$
Coefficients of linear combination \thetag{6.2} certainly depend on
$t_0$. For each particular $n$ they can be calculated explicitly.
These considerations lead to the following theorem.
\proclaim{Theorem 6.1} Suppose that Newtonian dynamical system on
$n$-dimensional Riemannian manifold $M$ is given. For any its
trajectory and for any choice of variation vector $\boldsymbol\tau$
on this trajectory corresponding function of deviation $\varphi(t)$
satisfies some linear homogeneous ordinary differential equation of
the order not greater than $2n$.
\endproclaim
\definition{Definition 6.1} We say that Newtonian dynamical system
on Riemannian manifold $M$ of the dimension $n\geqslant 2$ satisfies
{\bf weak normality} condition, if for each its trajectory there exists
some ordinary differential equation
$$
\ddot\varphi=\Cal A(t)\,\dot\varphi+\Cal B(t)\,\varphi\hskip -2em
\tag6.3
$$
such that any function of deviation $\varphi(t)$ corresponding to any
choice of variation vector $\boldsymbol\tau$ on that trajectory is
the solution of this equation.
\enddefinition
    Definition~6.1 is one of the main definitions in the theory of
dynamical systems admitting the normal shift. The definition of this
kind was first stated in \cite{38} (see also preprint \cite{34}).
It was found that weak normality condition leads to the system of
partial differential equations for the force field $\bold F$ of
Newtonian dynamical system. Derivation of these equations is our
main purpose in this section.\par
     Thus suppose, that Newtonian dynamical system with force field
$\bold F$ satisfies weak normality condition. Let's fix some trajectory
$p(t)$ and write corresponding differential equation \thetag{6.2} in
the following form:
$$
\ddot\varphi-\Cal A\,\dot\varphi-\Cal B\,\varphi=0.\hskip -2em
\tag6.4
$$
Let $p=p(t_0)$ be a point on the fixed trajectory, where velocity
vector doesn't vanish, i\.~e\. $\bold v(t_0)\neq 0$. In left hand side
of \thetag{6.3} we substitute the expression \thetag{4.2} for the
function of deviation $\varphi$, and substitute  the above expressions
\thetag{5.12} and \thetag{5.14} for its derivatives. After some
reductions we get
$$
\gathered
\left(\frac{\boldsymbol\alpha(\bold N)}{v}-\Cal A\right)\,
(\bold v\,|\,\nabla_t\boldsymbol\tau)+
\boldsymbol\alpha(\bold P\nabla_t\boldsymbol\tau)
+\boldsymbol\beta(\bold P\boldsymbol\tau)\,-\\
\vspace{1ex}
-\,\Cal A\,(\bold F\,|\,\bold P\,\boldsymbol\tau)
+\left(\frac{\boldsymbol\beta(\bold N)}{v}-\Cal A\,
\frac{(\bold F\,|\,\bold N)}{v}-\Cal B\right)\,
(\bold v\,|\,\boldsymbol\tau)=0.
\endgathered\hskip -2em
\tag6.5
$$
Let's consider the value of left hand side of \thetag{6.5} for
$t=t_0$. After fixing $t=t_0$ we have only one arbitrariness rest,
it consists in the choice of solution of the equation \thetag{3.11}.
This choice determines the value of $\boldsymbol\tau(t_0)$ and the
value of its derivative $\nabla_t\boldsymbol\tau(t_0)$. Left hand
side of \thetag{6.5} is linear in $\boldsymbol\tau$ and
$\nabla_t\boldsymbol\tau$. Therefore, considered for some fixed
$t=t_0$, \thetag{6.5} is an equality of linear functionals in
$\goth T^*$. The space of solutions of the equations \thetag{3.11}
on a given trajectory $p(t)$ is isomorphic to direct sum of two
copies of tangent space $T_p(M)$ at the point $p=p(t_0)$:
$$
\goth T\cong T_p(M)\oplus T_p(M).\hskip -2em
\tag6.6
$$
This is true, since $\goth T$ is parameterized by components of two
vectors $\boldsymbol\tau(t_0)\in T_p(M)$ and $\nabla_t\boldsymbol
\tau(t_0)\in T_p(M)$. Velocity vector $\bold v(t_0)\in T_p(M)$
subdivides each summand in formula \thetag{6.6} into direct sum of two
subspaces:
$$
T_p(M)=\langle\bold v\rangle\oplus\langle\bold v\rangle_{\sssize
\perp}.\hskip -2em
\tag6.7
$$
Here subspace $\langle\bold v\rangle$ is a linear span of velocity
vector $\bold v(t_0)\neq 0$, and $\langle\bold v\rangle_{\sssize\perp}$
is an orthogonal complement of $\langle\bold v\rangle$ to a whole
space $T_p(M)$ in the metric $g$. From \thetag{6.6} and \thetag{6.7}
we derive the expansion of $\goth T$ into a direct sum of four subspaces:
$$
\goth T\cong \langle\bold v\rangle\oplus\langle\bold v\rangle_{\sssize
\perp}\oplus\langle\bold v\rangle\oplus\langle\bold v\rangle_{\sssize
\perp}.\hskip -2em
\tag6.8
$$
The expansion \thetag{6.8} generates dual expansion in dual space
$\goth T^*$:
$$
\goth T^*\cong \langle\bold v\rangle^*\oplus\langle\bold v\rangle_{\sssize
\perp}^*\oplus\langle\bold v\rangle^*\oplus\langle\bold v\rangle_{\sssize
\perp}^*.\hskip -2em
\tag6.9
$$
Functionals $(\bold v\,|\,\nabla_t\boldsymbol\tau)$, \ $\boldsymbol
\alpha(\bold P\nabla_t\boldsymbol\tau)$, \ $\boldsymbol\beta(\bold P
\boldsymbol\tau)-\Cal A\,(\bold F\,|\,\bold P\boldsymbol\tau)$ and
$(\bold v\,|\,\boldsymbol\tau)$ in \thetag{6.5} belong to four
different subspaces in the expansion \thetag{6.9}. Therefore the
equation \thetag{6.5} is reduced to four separate equations. Two
of them determine $\Cal A$ and $\Cal B$ in \thetag{6.4}:
$$
\xalignat 2
&\quad\Cal A=\frac{\boldsymbol\alpha(\bold N)}{v},
&&\Cal B=\frac{\boldsymbol\beta(\bold N)}{v}-\Cal A\,
\frac{(\bold F\,|\,\bold N)}{v}.\hskip -3em
\tag6.10
\endxalignat
$$
Second pair of equation derived from \thetag{6.5} has the
following form:
$$
\xalignat 2
&\quad\boldsymbol\alpha(\bold P\nabla_t\boldsymbol\tau)=0,
&&\boldsymbol\beta(\bold P\boldsymbol\tau)=\Cal A\,(\bold F\,|\,
\bold P\,\boldsymbol\tau).\hskip -3em
\tag6.11
\endxalignat
$$
Let's calculate coefficients $\Cal A$ and $\Cal B$ in the equation
\thetag{6.3} by means of above formulas \thetag{6.10} and by means of
formulas \thetag{5.13} for the components of convectors
$\boldsymbol\alpha$ and $\boldsymbol\beta$:
$$
\align
&\Cal A=\sum^n_{i=1} \frac{2\,F_i\,N^i}{v}+
\sum^n_{i=1}\sum^n_{j=1} N^i\,N^j\,\tilde\nabla_iF_j,
\tag6.12\\
\vspace{1ex}
&\aligned
\Cal B&=\sum^n_{i=1}\sum^n_{j=1}\left(N^j\left(\nabla_jF_i+
\nabla_iF_j\right)+\frac{F^j\,\tilde\nabla_jF_i}{v}\right)N^i\,-\\
&-\,\sum^n_{i=1}\sum^n_{k=1}\left(\frac{2\,F_i\,N^i}{v}
+\shave{\sum^n_{j=1}}N^i\,N^j\,\tilde\nabla_iF_j\right)
\frac{F_k\,N^k}{v}.
\endaligned\hskip -3em
\tag6.13
\endalign
$$\par
    Now consider the equations \thetag{6.11}. Using \thetag{5.13},
one can write them in local coordinates. Afterwards one should take
into account arbitrariness of vectors $\boldsymbol\tau$ and
$\nabla_t\boldsymbol\tau$. As a result first equation \thetag{6.11}
is written in the following form:
$$
\sum^n_{i=1}\left(2\,F_i+\shave{\sum^n_{j=1}}v^j\,
\tilde\nabla_iF_j\right)P^i_k=0.\hskip -2em
\tag6.14
$$
In order to transform another equation \thetag{6.11} we use formula
\thetag{6.12} for the coefficient $\Cal A$ in the equation \thetag{6.3}.
This yields
$$
\aligned
\sum^n_{i=1}\left(\shave{\sum^n_{j=1}}v^j\left(\nabla_jF_i+
\nabla_iF_j\right)+\shave{\sum^n_{j=1}} F^j\,\tilde\nabla_jF_i\right)
P^i_k=\\
\vspace{1ex}
=\sum^n_{r=1}\left(\shave{\sum^n_{i=1}}\frac{2\,F_i\,N^i}{v}+
\shave{\sum^n_{i=1}}\shave{\sum^n_{j=1}}N^i\,N^j\,\tilde\nabla_iF_j
\right)F_r\,P^r_k.
\endaligned\hskip -3em
\tag6.15
$$
All done. We are only to do some slight cosmetic transformations in the
equation \thetag{6.14} and \thetag{6.15}. The equation \thetag{6.14} is
written as
$$
\sum^n_{i=1}\left(v^{-1}\,F_i+\shave{\sum^n_{j=1}}
\tilde\nabla_i\left(N^j\,F_j\right)\right)P^i_k=0.\hskip -2em
\tag6.16
$$
In deriving \thetag{6.16} we took into account $v^j=v\,N^j$, and used
second relationship \thetag{5.7} together with well-known property
$\bold P^2=\bold P$ of projector $\bold P$. In equation \thetag{6.15}
we shall only transpose some terms:
$$
\pagebreak
\aligned
&\sum^n_{i=1}\sum^n_{j=1}\left(\nabla_iF_j+\nabla_jF_i-2\,v^{-1}
\,F_i\,F_j\right)N^j\,P^i_k\,+\\
\vspace{1ex}
&+\sum^n_{i=1}\sum^n_{j=1}\left(\frac{F^j\,\tilde\nabla_jF_i}{v}
-\sum^n_{r=1}\frac{N^r\,N^j\,\tilde\nabla_jF_r}{v}\,F_i\right)P^i_k=0.
\endaligned\hskip -3em
\tag6.17
$$
The equations \thetag{6.16} and \thetag{6.17}, which we have derived just
now, are called {\bf weak normality} equations. In this form they were
first obtained in \cite{58} (see also \cite{75}). The above derivation of
these equations proves the following theorem.
\proclaim{Theorem 6.2} Newtonian dynamical system on Riemannian manifold
of the dimension $n\geqslant 2$ satisfies weak normality condition if and
only if its force field satisfies the differential equations \thetag{6.16}
and \thetag{6.17} at all points $q=(p,\bold v)$ of tangent bundle $TM$,
except for those, where $\bold v=0$.
\endproclaim
    The equations of weak normality \thetag{6.16} and \thetag{6.17}
form the system of partial differential equations for the force field
$\bold F$. The number of these equations is $2n$. However, projector
matrix $P^i_k$ is degenerate, its rank is $n-1$. Therefore actual number
of independent equations in \thetag{6.16} and \thetag{6.17} is $2n-2$.
The problem of their compatibility is studied in next two Chapters, where
we construct great many of their solutions. Here we note only that
identically zero force field $\bold F$ is the solution of the equations
\thetag{6.16} and \thetag{6.17}. This proves the following theorem.
\proclaim{Theorem 6.3} Geodesic flow on Riemannian manifold $M$ of the
dimension $n\geqslant 2$ is a Newtonian dynamical system satisfying weak
normality condition.
\endproclaim
\head
\S\,7. Cauchy problem for the function of deviation.
\endhead
    Suppose that $M$ is a Riemannian manifold with Newtonian dynamical
system satisfying weak normality condition on it. Let's use this system
to define the shift $f_t\!:\,S\to S_t$ for some hypersurface $S$ from
transformation class on $M$ (see definition~2.1). Functions of deviation
on the trajectories of shift satisfy the differential equations
\thetag{6.3}. These are linear homogeneous differential equations of the
second order:
$$
\ddot\varphi=\Cal A(t)\,\dot\varphi+\Cal B(t)\,\varphi.\hskip -2em
\tag7.1
$$
Coefficients $\Cal A(t)$ and $\Cal B(t)$ can be different for different
trajectories of shift. Upon choosing local coordinates $u^1,\,\ldots,\,
u^{n-1}$ on $S$ and local coordinates $x^1,\,\ldots,\,x^n$ in $M$, we can
define shift functions \thetag{3.1}. Their derivatives determine vectors
of variation \thetag{3.3}, which coincide with coordinate tangent vectors
to hypersurfaces $S_t$:
$$
\boldsymbol\tau_i=\sum^n_{j=1}\frac{\partial x^j}{\partial u^i}\,
\frac{\partial}{\partial x^j}\text{, \ where \ }i=1,\,\ldots,\,n-1.
\hskip -2em
\tag7.2
$$
By differentiating shift functions \thetag{3.1} with respect to time
variable $t$ for the fixed values of  $u^1,\,\ldots,\,u^{n-1}$ we
obtain the velocity vector on the trajectories of shift:
$$
\bold v=\sum^n_{j=1}\dot x^j\,\frac{\partial}{\partial x^j}.
\hskip -2em
\tag7.3
$$
Scalar products of velocity vector $\bold v$ from \thetag{7.3}
with variation vectors $\boldsymbol\tau_1,\,\ldots,\,\boldsymbol
\tau_{n-1}$ from \thetag{7.2} define the set functions of deviations
$$
\varphi_1,\ \ldots,\ \varphi_{n-1}.\hskip -2em
\tag7.4
$$
This is in complete agreement with formulas \thetag{4.1}. Each function
\thetag{7.4} satisfies the equations \thetag{7.1} for fixed values of
variables $u^1,\,\ldots,\,u^{n-1}$. Normality condition for the shift
$f_t\!:\,S\to S_t$ consists in identical vanishing of all functions
\thetag{7.4} (see definition~2.3). In order to satisfy this condition
it is sufficient to provide vanishing of functions $\varphi_k$ and their
first derivatives $\dot\varphi_k$ at the initial instant of time $t=0$,
i\.~e\. at the points of initial hypersurface $S$. Thus, we are to set up
the following Cauchy problem for differential equations \thetag{7.1}
on trajectories of shift:
$$
\xalignat 2
&\quad\varphi_k(t)\,\hbox{\vrule height 8pt depth 8pt width
0.5pt}_{\,t=0}=0,
&&\dot\varphi_k(t)\,\hbox{\vrule height 8pt depth 8pt width
0.5pt}_{\,t=0}=0.\hskip -2em
\tag7.5
\endxalignat
$$
First part of initial conditions \thetag{7.5} appears to be fulfilled
unconditionally due to the initial data \thetag{2.2} that define
trajectories of shift $f_t\!:\,S\to S_t$. Indeed, according to
\thetag{2.2} vector of initial velocity $\bold v(0)$ is directed
along the normal vector of $S$. Therefore it is perpendicular to the
vectors $\boldsymbol\tau_1(0),\,\ldots,\,\boldsymbol\tau_{n-1}(0)$,
which are tangent to $S$.\par
    Let's study in more details the second pert of initial conditions
\thetag{7.5}. In order to do it we use formula \thetag{4.7}. Here this
formula is written as
$$
\dot\varphi_k=(\bold F\,|\,\boldsymbol\tau_k)+(\bold v\,|\,
\nabla_t\boldsymbol\tau_k).\hskip -2em
\tag7.6
$$
In order to calculate vector  $\nabla_t\boldsymbol\tau_k$ in \thetag{7.6}
we apply formula \thetag{5.13} from Chapter
\uppercase\expandafter{\romannumeral 4} and take unto account that here
$u=u^k$. Tensor field of torsion $\bold T$ for Riemannian connection
\thetag{1.2} is equal to zero, therefore
$$
\nabla_t\boldsymbol\tau_k=\nabla_{u^k}\bold v.\hskip -2em
\tag7.7
$$
Let's write the equality \thetag{7.7} for $t=0$, taking into
account initial data \thetag{2.2}. In order to calculate covariant
derivative $\nabla_{u^k}\bold u$ in right hand side of \thetag{7.7}
one should only know the values of velocity vector on initial
hypersurface $S$. This is clear, for instance, due to formula
\thetag{5.8} from Chapter \uppercase\expandafter{\romannumeral 4}.
Then
$$
\nabla_t\boldsymbol\tau_k\,\hbox{\vrule height 8pt depth 8pt width
0.5pt}_{\,t=0}=\nabla_{u^k}(\nu\cdot\bold n)=\frac{\partial\nu}
{\partial u^k}\cdot\bold n+\nu\cdot\nabla_{u^k}\bold n.\hskip -2em
\tag7.8
$$
Covariant derivatives $\nabla_{u^k}\bold n$ are natural in the
theory of hypersurfaces embedded into Riemannian manifold. They
are encountered in derivational formulas of Weingarten (see \cite{16},
\cite{32}, \cite{84}, or \cite{77}). In our case for the hypersurface
$S$ in $M$ derivational formulas of Weingarten are written as follows:
$$
\align
\nabla_{u^k}\boldsymbol\tau_r&=\sum^{n-1}_{m=1}\theta^m_{kr}\,
\boldsymbol\tau_m+b_{kr}\,\bold n,\hskip -3em
\tag7.9\\
\nabla_{u^k}\bold n&=-\sum^{n-1}_{m=1}b^m_k\,\boldsymbol\tau_m.
\hskip -3em
\tag7.10
\endalign
$$
If we denote by $\rho_{ij}$ \pagebreak components of metric tensor 
$\boldsymbol\rho$
for induced Riemannian metric $\boldsymbol\rho$ in local coordinates $u^1,
\,\ldots,\,u^{n-1}$ on $S$, then quantities $\theta^m_{kr}$ in \thetag{7.9}
are the components of corresponding metric connection for the metric
$\boldsymbol\rho$. They are calculated by formula analogous to formula
\thetag{1.2} for $\Gamma^k_{ij}$:
$$
\theta^k_{ij}=\frac{1}{2}\sum^{n-1}_{s=1}\rho^{ks}\left(\frac{\partial
\rho_{sj}}{\partial u^i}+\frac{\partial\rho_{is}}{\partial u^j}-\frac{
\partial\rho_{ij}}{\partial u^s}\right).\hskip -3em
\tag7.11
$$
Quantities $b_{kr}$ in \thetag{7.9} are the components of inner tensor
field $\bold b$ of type $(0,2)$ in $S$ which is called {\bf second
fundamental form} of hypersurface $S$. They are symmetric, i\.~e\.
they do not change in transposing indices:
$$
b_{kr}=b_{rk}.\hskip -3em
\tag7.12
$$
Quantities $b^m_k$ in formula \thetag{7.10} are derived from $b_{kr}$
by index raising procedure:
$$
b^m_k=\sum^{n-1}_{r=1}b_{kr}\,\rho^{rm}.\hskip -3em
\tag7.13
$$
They are the components of inner tensor field of type $(1,1)$ in
submanifold $S$. Such field is  an  operator  version  of  second 
fundamental
form $\bold b$ in $S$.\par
     Let's apply derivational formula of Weingarten \thetag{7.10} in
order to express covariant derivative $\nabla_{u^k}\bold n$ in formula
\thetag{7.8}. This yields
$$
\nabla_t\boldsymbol\tau_k\,\hbox{\vrule height 8pt depth 8pt width
0.5pt}_{\,t=0}=\frac{\partial\nu}{\partial u^k}\cdot\bold n-\nu\cdot
\sum^{n-1}_{m=1}b^m_k\,\boldsymbol\tau_m.\hskip -3em
\tag7.14
$$
Now we substitute \thetag{7.14} into \thetag{7.6} and calculate
derivative $\dot\varphi_k$ for initial instant of time $t=t_0$
on initial hypersurface $S$:
$$
\dot\varphi_k\,\hbox{\vrule height 8pt depth 8pt width 0.5pt}_{\,t=0}
=(\bold F\,|\,\boldsymbol\tau_k)+\nu\,\frac{\partial\nu}{\partial u^k}.
\hskip -3em
\tag7.15
$$
In deriving \thetag{7.15} we took into account second part of the
relationships \thetag{2.2} that determine initial value of velocity
vector on $S$.\par
    Let's substitute the above formula \thetag{7.15} expressing initial
value of derivative $\dot\varphi$ into the equations \thetag{7.5}. This
yields the equations
$$
\frac{\partial\nu}{\partial u^k}=-\nu^{-1}\,
(\bold F\,|\,\boldsymbol\tau_k)
\hskip -3em
\tag7.16
$$
for the function $\nu(u^1,\ldots,u^{n-1})$. If $\nu(u^1,\ldots,u^{n-1})$
is the solution of equations \thetag{7.16}, then for all function of
deviation \thetag{7.4} initial conditions \thetag{7.5} are fulfilled.
In the case of dynamical system satisfying weak normality condition
(see definition~6.1) it is sufficient to assure the normality of
shift $f_t\!:\,S\to S_t$.
\head
\S\,8. Pfaff equations for the modulus of initial velocity.
\endhead
     The equations \thetag{7.16} determine two essential different
cases: {\bf two dimensional case}, when $\dim M=n=2$, when we have
only one equation \thetag{7.16}, and {\bf multidimensional case},
when $\dim M=n\geqslant 3$, when the number of equations is greater
than $1$. We intentionally omit here two dimensional case, since it
is considered in {\bf other thesis} \cite{36}.\par
     Let's consider the equations \thetag{7.16} in multidimensional
case. All quantities in right hand side of thetag{7.16} are localized
on initial hypersurface $S$. Force vector $\bold F$ depends on the
point of tangent bundle $q=(p,\bold v)$. Substituting $\bold F$ into
\thetag{7.16} we choose $p\in S$ and $\bold v=\nu\cdot\bold n(p)$ in
accordance with initial data \thetag{2.2}. In local coordinates
$u^1,\,\ldots,\,u^{n-1}$ the point $p\in S$ and coordinate tangent
vectors  $\boldsymbol\tau_k$ are parameterized by these coordinates:
$p=p(u^1,\ldots,u^{n-1})$ and $\boldsymbol\tau_k=\boldsymbol\tau_k(u^1,
\ldots,u^{n-1})$. Therefore the equations \thetag{7.16} has the
following structure:
$$
\frac{\partial\nu}{\partial u^k}=\psi_k(\nu,u^1,\ldots,u^{n-1})
\text{, \ where \ }k=1,\,\ldots,\,n-1.
\hskip -3em
\tag8.1
$$
For $n\geqslant 3$ the equations \thetag{8.1} form {\bf complete
system of Pfaff equations}. Central point of the theory of Pfaff
equations is the {\bf problem of compatibility}. Suppose that
$\nu=\nu(u^1,\ldots,u^{n-1})$ is a solution of the equations
\thetag{8.1}. Let's calculate second order partial derivative
$\nu$ with respect to variables $u^k$ and $u^r$ by virtue of
the equations \thetag{8.1}. This can be done in two ways:
$$
\align
\frac{\partial^2\nu}{\partial u^k\,\partial u^r}&=
\frac{\partial\psi_k}{\partial u^r}+\frac{\partial\psi_k}{\partial\nu}
\,\psi_r=\vartheta_{kr}(\nu,u^1,\ldots,u^{n-1}),\hskip -2em
\tag8.2\\
\vspace{1ex}
\frac{\partial^2\nu}{\partial u^r\,\partial u^k}&=
\frac{\partial\psi_r}{\partial u^k}+\frac{\partial\psi_r}{\partial\nu}
\,\psi_k=\vartheta_{rk}(\nu,u^1,\ldots,u^{n-1}).\hskip -2em
\tag8.3
\endalign
$$
Functions $\vartheta_{kr}(\nu,u^1,\ldots,u^{n-1})$ and $\vartheta_{kr}
(\nu,u^1,\ldots,u^{n-1})$ in right hand sides of \thetag{8.2} and \thetag{8.3}
can be calculated regardless to whether has the system of Pfaff equations
thetag{8.1} some solution or it hasn't. In the case, when it has, the
derivatives \thetag{8.2} and \thetag{8.3} are equal to each other. Equating
\thetag{8.2} and \thetag{8.3}, we get
$$
\aligned
&\vartheta_{kr}(\nu(u^1,\ldots,u^{n-1}),u^1,\ldots,u^{n-1})=\\
\vspace{1ex}
&=\vartheta_{rk}(\nu(u^1,\ldots,u^{n-1}),u^1,\ldots,u^{n-1}).
\endaligned\hskip -2em
\tag8.4
$$
In general, the equality \thetag{8.4} doesn't  provide  identical 
coincidence
of the functions $\vartheta_{kr}$ and $\vartheta_{rk}$. In order to test
whether it is fulfilled or not we are to have some solution $\nu(u^1,
\ldots,u^{n-1})$ of the equations \thetag{8.1}. But there is some situation
when we can escape this difficulty.
\definition{Definition 8.1} System of Pfaff equations \thetag{8.1} is
called {\bf compatible}, if for any fixed point $(\nu_0,\,u^1_0,\,\ldots,
\,u^{n-1}_0)$ from the domain of the functions $\psi_k$ \pagebreak there
is a solution $\nu=\nu(u^1,\ldots,u^{n-1})$ of the equations \thetag{8.1}
defined in some neighborhood of the point $(u^1_0,\ldots,u^{n-1}_0)$ and
normalized by the condition 
$$
\nu(u^1_0,\ldots,u^{n-1}_0)=\nu_0.\hskip -2em
\tag8.5
$$
\enddefinition
     If the system of Pfaff equations \thetag{8.1} is compatible, then
from \thetag{8.5} and \thetag{8.4} we get that functions $\vartheta_{kr}$
and $\vartheta_{rk}$ do identically coincide
$$
\vartheta_{kr}(\nu,u^1,\ldots,u^{n-1})=\vartheta_{rk}(\nu,u^1,\ldots,
u^{n-1}).\hskip -2em
\tag8.6
$$
The relationships \thetag{8.6} are called {\bf compatibility conditions}
for complete system of Pfaff equations \thetag{8.1}. Let's write them
in more explicit form:
$$
\frac{\partial\psi_k}{\partial u^r}+\frac{\partial\psi_k}{\partial\nu}
\,\psi_r=\frac{\partial\psi_r}{\partial u^k}+\frac{\partial\psi_r}
{\partial\nu}\,\psi_k.\hskip -2em
\tag8.7
$$
\proclaim{Theorem 8.1} Compatibility conditions expressed by relationships
\thetag{8.7}, where $k=1,\,\ldots,\,n-1$ and $r=1,\,\ldots,\,n-1$, are
necessary and sufficient for the system of Pfaff equations \thetag{8.1} to
be compatible in the sense of definition~8.1.
\endproclaim
\demo{Proof} The necessity was proved by derivations of relationships
\thetag{8.7}. We are to prove sufficiency. Let's do it by induction in
$n$. For $n=2$ we have one equation \thetag{8.1} which is ordinary
differential equation with independent variable $u^1$. Normalization
\thetag{8.5} is a Cauchy problem for this equation. Compatibility
conditions \thetag{8.7} in this case are identically fulfilled, since
$k=r=1$. Thus, for $n=2$ the proposition of theorem follows from
local solvability of Cauchy problem for ordinary differential equation
(see ``{existence and uniqueness}'' theorem in \cite{82} or \cite{83}).
This is the base of induction.\par
    Let's do the step of induction from $n-1$ to $n$. If we exclude the
last equation from the system \thetag{8.1} and if we fix the last
variable $u^{n-1}=u^{n-1}_0$, then we get into the situation when
inductive hypothesis is applicable. This yields the function
$$
\nu=\nu(u^1,\ldots,u^{n-2})\hskip -2em
\tag8.8
$$
that satisfies the system of first $(n-2)$ Pfaff equations in
\thetag{8.1}, and that satisfies normalizing condition \thetag{8.5}
on the hyperplane $u^{n-1}=u^{n-1}_0$. Now let's consider again the
last equation in \thetag{8.1}, treating it as an ordinary differential
equation with independent variable $u^{n-1}$, and let's set up the
Cauchy problem
$$
\nu\,\hbox{\vrule height 8pt depth 8pt width 0.5pt}_{\,u^{n-1}=
u^{n-1}_0}=\nu(u^1,\ldots,u^{n-2}),\hskip -2em
\tag8.9
$$
using the function \thetag{8.8} as initial data. We get the required
function $\nu(u^1,\ldots,u^{n-1})$ as a solution of the above Cauchy
problem with initial data \thetag{8.9}.\par
    Function $\nu(u^1,\ldots,u^{n-1})$ satisfies the last equation
in the system \thetag{8.1} by construction. Moreover, it satisfies
other equations \thetag{8.1} on the hyperplane $u^{n-1}=u^{n-1}_0$
by inductive hypothesis. Suppose that  $1\leqslant k\leqslant n-2$.
Let's show that $k$-th equation \thetag{8.1} is fulfilled for
all values of arguments in $\nu(u^1,\ldots,u^{n-1})$ where this
function is defined. Consider the following functions:
$$
\align
x_k(u^1,\ldots,u^{n-1})&=\frac{\partial\nu}{\partial u^k},\hskip -3em
\tag8.10\\
\vspace{1.5ex}
y_k(u^1,\ldots,u^{n-1})&=\psi_k(\nu,u^1,\ldots,u^{n-1}).\hskip -3em
\tag8.11
\endalign
$$
Let's calculate the derivatives of \thetag{8.10} and \thetag{8.11}
with respect to the variable $u^{n-1}$. For the first of these two
functions we have 
$$
\frac{\partial x_k}{\partial u^{n-1}}=\frac{\partial^2\nu}
{\partial u^k\,\partial u^{n-1}}=\frac{\partial}{\partial u^k}
\left(\psi_{n-1}(\nu,u^1,\ldots,u^{n-1})\right).
$$
By differentiating composite function in the above expression we get
$$
\frac{\partial x_k}{\partial u^{n-1}}=\frac{\partial\psi_{n-1}}
{\partial u^k}+\frac{\partial\psi_{n-1}}{\partial\nu}\,x^k.
\hskip -3em
\tag8.12
$$
By means of analogous differentiation for the function \thetag{8.11} we
obtain
$$
\frac{\partial y_k}{\partial u^{n-1}}=\frac{\partial\psi_k}
{\partial u^{n-1}}+\frac{\partial\psi_k}{\partial\nu}\,\psi_{n-1}.
\hskip -3em
\tag8.13
$$
Let's transform \thetag{8.13} with the help of compatibility condition
\thetag{8.7}, taking $r=n-1$ in it. This yields the relationship
$$
\frac{\partial y_k}{\partial u^{n-1}}=\frac{\partial\psi_{n-1}}
{\partial u^k}+\frac{\partial\psi_{n-1}}{\partial\nu}\,y^k.
\hskip -3em
\tag8.14
$$
Comparing \thetag{8.12} and \thetag{8.14}, we see that functions
$x_k$ and $y_k$ satisfy identical linear ordinary differential equations
of the first order with to independent variable $u^{n-1}$. For
$u^{n-1}=u^{n-1}_0$ their values are equal to each other by
inductive hypothesis. Therefore $x_k=y_k$ for all $u^1,\,\ldots,\,u^{n-1}$.
This means that function $\nu(u^1,\ldots,u^{n-1})$ constructed above
is the solution of whole system of Pfaff equations \thetag{8.1}. And
it satisfies normalizing condition \thetag{8.5} as well. Theorem is
proved.\qed\enddemo
    Let's apply the theorem~8.1 to the equations \thetag{7.16}, and
let's find compatibility condition for them. First calculate partial
derivatives \thetag{8.2} and \thetag{8.3} for the function
$\nu(u^1,\ldots,u^{n-1})$ by virtue of the equations \thetag{7.16}:
$$
\align
&\aligned
&\frac{\partial^2\nu}{\partial u^k\,\partial u^r}=-\nabla_{u^r}
\left(\nu^{-1}\,(\bold F\,|\,\boldsymbol\tau_k)\right)
=\frac{(\bold F\,|\,\boldsymbol\tau_k)}{\nu^2}\,\frac{\partial\nu}
{\partial u^r}-\frac{(\nabla_{u^r}\bold F\,|\,\boldsymbol\tau_k)}{v}\,-\\
\vspace{2ex}
&-\,\frac{(\bold F\,|\,\nabla_{u^r}\boldsymbol\tau_k)}{v}=
-\frac{(\bold F\,|\,\boldsymbol\tau_k)\,(\bold F\,|\,
\boldsymbol\tau_r)}{\nu^3}-\frac{(\nabla_{u^r}\bold F\,|\,
\boldsymbol\tau_k)}{\nu}-\frac{(\bold F\,|\,\nabla_{u^r}
\boldsymbol\tau_k)}{\nu},
\endaligned\hskip -3em
\tag8.15\\
\displaybreak
&\aligned
&\frac{\partial^2\nu}{\partial u^r\,\partial u^k}=-\nabla_{u^k}
\left(\nu^{-1}\,(\bold F\,|\,\boldsymbol\tau_r)\right)
=\frac{(\bold F\,|\,\boldsymbol\tau_r)}{\nu^2}\,\frac{\partial\nu}
{\partial u^k}-\frac{(\nabla_{u^k}\bold F\,|\,\boldsymbol\tau_r)}{v}\,-\\
\vspace{2ex}
&-\frac{(\bold F\,|\,\nabla_{u^k}\boldsymbol\tau_r)}{v}=
-\frac{(\bold F\,|\,\boldsymbol\tau_r)\,(\bold F\,|\,
\boldsymbol\tau_k)}{\nu^3}-\frac{(\nabla_{u^k}\bold F\,|\,
\boldsymbol\tau_r)}{\nu}-\frac{(\bold F\,|\,\nabla_{u^k}
\boldsymbol\tau_r)}{\nu}.
\endaligned\hskip -3em
\tag8.16
\endalign
$$
Let's subtract \thetag{8.16} from \thetag{8.15}, taking into account
that $\nabla_{u^r}\boldsymbol\tau_k=\nabla_{u^k}\boldsymbol\tau_r$.
The latter equality is the consequence of derivational formulas
\thetag{7.9} due to symmetry of components of metric connection
\thetag{7.11}, and due to the relationship \thetag{7.12} expressing
symmetry of second fundamental form. Taking into account all above
facts, we get
$$
(\nabla_{u^r}\bold F\,|\,\boldsymbol\tau_k)=(\nabla_{u^k}\bold F\,|\,
\boldsymbol\tau_r).\hskip -3em
\tag8.17
$$
In order to calculate covariant derivative $\nabla_{u^r}\bold F$
in the equality \thetag{8.17} we apply the formula \thetag{5.19}
from Chapter \uppercase\expandafter{\romannumeral 4}:
$$
\nabla_{u^r}\bold F=C(\nabla_{u^r}\bold v\otimes\tilde\nabla\bold F)
+C(\boldsymbol\tau_r\otimes\nabla\bold F).\hskip -3em
\tag8.18
$$
In calculating $\nabla_{u^r}\bold v$ we take into account the
equality $\bold v=\nu\cdot\bold n$ for the initial velocity on
the initial hypersurface $S$. Then we get 
$$
\nabla_{u^r}\bold v=\nabla_{u^r}(\nu\cdot\bold n)=-\frac{(\bold F\,|\,
\boldsymbol\tau_r)}{\nu}\cdot\bold n-\nu\cdot\sum^{n-1}_{m=1}b^m_r\,
\boldsymbol\tau_m.\hskip -3em
\tag8.19
$$
In obtaining \thetag{8.19} we used derivational formulas \thetag{7.10}
and the equation \thetag{7.16} for $\nu$ by itself. Now let's substitute
\thetag{8.19} into \thetag{8.18}, and let's use the resulting expression
for to transform the equation \thetag{8.17}:
$$
\gather
\frac{(\bold F\,|\,\boldsymbol\tau_k)(C(\bold n\otimes\tilde\nabla\bold F)
\,|\,\boldsymbol\tau_r)}{\nu}-\frac{(\bold F\,|\,\boldsymbol\tau_r)(C(\bold
n\otimes\tilde\nabla\bold F)\,|\,\boldsymbol\tau_k)}{\nu}\,+\hskip -3em\\
\vspace{2ex}
+\,(C(\boldsymbol\tau_r\otimes\nabla\bold F)\,|\,\boldsymbol\tau_k)
-(C(\boldsymbol\tau_k\otimes\nabla\bold F)\,|\,\boldsymbol\tau_r)=
\hskip -3em
\tag8.20\\
\vspace{1ex}
=\sum^{n-1}_{m=1}\nu\left(b^m_r\,(C(\boldsymbol\tau_m\otimes\tilde\nabla
\bold F)\,|\,\boldsymbol\tau_k)-b^m_k\,(C(\boldsymbol\tau_m\otimes\tilde
\nabla\bold F)\,|\,\boldsymbol\tau_r)\right).\hskip -3em
\endgather
$$
The relationship \thetag{8.20} is the compatibility condition \thetag{8.7}
applied to the Pfaff equations \thetag{7.16}. Let's write it in local
coordinates on $M$:
$$
\gathered
\sum^n_{i=1}\sum^n_{j=1}\tau^i_k\,\tau^j_r\left(\,\shave{\sum^n_{m=1}}
n^m\,\frac{F_i\,\tilde\nabla_mF_j-F_j\,\tilde\nabla_mF_i}{\nu}+
\nabla_jF_i-\nabla_iF_j\right)=\\
\vspace{2ex}
=\sum^n_{i=1}\sum^n_{j=1}\sum^{n-1}_{m=1}\nu\left(b^m_r\,
\tau^j_m\tilde\nabla_jF_i\,\tau^i_k-b^m_k\,\tau^j_m\tilde\nabla_jF_i
\,\tau^i_r\right).
\endgathered\hskip -1.5em
\tag8.21
$$
Here $n^m$ are the components of unitary normal vector to $S$, and
$\tau^i_k$ are components of the vector $\boldsymbol\tau_k$ in local
map $(U,x^1,\ldots,x^n)$ on $M$:
$$
\boldsymbol\tau_k=\sum^n_{i=1}\tau^i_k\,
\frac{\partial}{\partial x^i}\text{, \ where \ }k=1,\,\ldots,\,n-1.
\hskip -3em
\tag8.22
$$
Variation vectors \thetag{8.22} tangent to $S$ and the components
$b^m_k$ of inner operator field $\bold b$ from \thetag{7.10} are
determined by the choice of local coordinates $u^1,\,\ldots,\,u^{n-1}$
on the hypersurface $S$. The equations \thetag{8.21} expressing
compatibility condition \thetag{8.20} in local coordinates appear
to be more convenient for the further analysis.
\head
\S\,9. The additional normality condition.
\endhead
    Normalization condition \thetag{8.5} in the theory of Pfaff
equations is similar to initial data in cauchy problem for ordinary
differential equations. Theorem~8.1 is analog of the theorem on
local solvability of Cauchy problem for ordinary differential equations.
The condition \thetag{8.5} prompt us normalization condition for
the function $\nu(p)$ that determines modulus of initial velocity in
\thetag{2.2}. Let $p_0$ be some fixed point on hypersurface $S$ in $M$,
and let $\nu_0\neq 0$ be some fixed real number. Then we set
$$
\nu(p)\,\hbox{\vrule height 8pt depth 8pt width 0.5pt}_{\,p=p_0}=
\nu_0\hskip -2em
\tag9.1
$$
This equality \thetag{9.1} is called {\bf normalization condition}
for the function $\nu(p)$ on $S$.
\definition{Definition 9.1} Suppose that Riemannian manifold $M$ of
the dimension $n\geqslant 2$ is equipped with some Newtonian dynamical
system, which is used to arrange the shift of hypersurfaces along its
trajectories according to the initial data \thetag{2.2}. We say that
this dynamical system satisfies {\bf additional normality} condition,
if for any hypersurface $S$ in $M$, for any point $p_0$ on $S$, and
for any real number $\nu_0\neq 0$ there exists some smaller piece $S'$
of hypersurface $S$ belonging to transformation class and containing
the point $p_0$, and there exists some function $\nu(p)$ from
transformation class on $S'$ normalized by the condition \thetag{9.1}
and such that, when substituted to \thetag{2.2}, it defines the
shift $f_t\!:\,S'\to S'_t$ such that 
$$
\dot\varphi_k(t)\,\hbox{\vrule height 8pt depth 8pt width
0.5pt}_{\,t=0}=0\text{, \ \ where \ \ \ }k=1,\,\ldots,\,n-1.
\hskip -2em
\tag9.2
$$
Here $\varphi_k=(\boldsymbol\tau_k\,|\,\bold v)$ are the functions of
deviation defined by some choice of local coordinates $u^1,\,\ldots,
\,u^{n-1}$ on $S$ and corresponding vectors of variation $\boldsymbol
\tau_1,\,\ldots,\,\boldsymbol\tau_{n-1}$.
\enddefinition
     Note that \thetag{9.2} is a second part of the equalities
\thetag{7.5}. If additional normality condition is fulfilled, then
equalities \thetag{7.5} are also fulfilled, since first part of
these equalities holds unconditionally due to the initial data
\thetag{2.2}.\par
     Functions of deviations $\varphi_k$ in \thetag{9.2} depend
on the choice of local coordinates $u^1,\,\ldots,\,u^{n-1}$ on
$S$. But additional normality condition from definition~9.1
is fulfilled or not fulfilled regardless to any particular choice
of local coordinates. Indeed, \pagebreak if we change local coordinates
$u^1,\,\ldots,\,u^{n-1}$ on $S$ for other local coordinates $\tilde u^1,
\,\ldots,\,\tilde u^{n-1}$, then corresponding vectors of variation
are transformed as follows
$$
\tilde{\boldsymbol\tau}_i=\sum^{n-1}_{k=1}\sigma^k_i\,\boldsymbol\tau_k,
\hskip -2em
\tag9.3
$$
Here $\sigma^k_i$ are components of transition matrix given by
$$
\sigma^k_i=\frac{\partial u^k}{\partial\tilde u^i}.
$$
They do not depend on $t$. The matter is that local coordinates $u^1,\,
\ldots,\,u^{n-1}$ and $\tilde u^1,\,\ldots,\,\tilde u^{n-1}$ are transferred
from $S'$ to $S'_t$ by bijective map of shift $f_t\!:\,S'\to S'_t$.
Transition functions 
$$
\xalignat 2
&\tilde u^1=\tilde u^1(u^1,\ldots,u^{n-1}), &&u^1=u^1(\tilde u^1,
\ldots,\tilde u^{n-1}),\\
&.\ .\ .\ .\ .\ .\ .\ .\ .\ .\ .\ .\ .\ .\ .\ .\ .\ .\ .\ .\ &&.\ .\ .\ .\
.\ .\ .\ .\ .\ .\ .\ .\ .\ .\ .\ .\ .\ .\ .\ .\ \\
&\tilde u^{n-1}=\tilde u^{n-1}(u^1,\ldots,u^{n-1}),\quad &&u^{n-1}=
u^{n-1}(\tilde u^1,\ldots,\tilde u^{n-1})
\endxalignat
$$
thereby remain same as on initial hypersurface $S'$. Transformation rule
\thetag{9.3} then is the direct consequence of formula \thetag{7.2}, which,
in turn, is the definition of variation vectors $\boldsymbol\tau_1,\,
\ldots,\,\boldsymbol\tau_{n-1}$.\par
    Let's substitute \thetag{9.3} into \thetag{4.1}. As a result we obtain
the transformation rule for the functions of deviation relating $\varphi_1,
\,\ldots,\,\varphi_{n-1}$ and $\tilde\varphi_1,\,\ldots,\,\tilde
\varphi_{n-1}$:
$$
\tilde\varphi_i=\sum^{n-1}_{k=1}\sigma^k_i\,\varphi_k.\hskip -2em
\tag9.4
$$
Differentiating \thetag{9.4} with respect to $t$ and taking into account
$\dot\sigma^k_i=0$, we get
$$
\DotTildePhi_i=\sum^{n-1}_{k=1}\sigma^k_i\,\dot\varphi_k.
\hskip -2em
\tag9.5
$$
The relationships \thetag{9.5} shows that, once fulfilled for some
particular choice of local coordinates on $S$, the relationships
\thetag{9.2} will remain true for any other choice of local
coordinates on $S$.\par
     The statement of additional normality condition in definition~9.1
exploits some hypersurface $S$ and some local coordinates on it. But,
as we have seen above, these objects are auxiliary tools only.
Additional normality condition, when it is fulfilled, is a property
of dynamical system by itself. Similar to weak normality condition
from definition~6.1, it can be written in form of partial differential
equations for the force field $\bold F$ of Newtonian dynamical system.
Exception is two-dimensional case $\dim M=n=2$.
\proclaim{Theorem 9.1}In two-dimensional case $\dim M=2$ additional
normality condition is fulfilled unconditionally for arbitrary
\pagebreak Newtonian dynamical system on $M$.
\endproclaim
    We shall not prove and comment this theorem, since two-dimensional
case is considered in {\bf other thesis} \cite{36}.\par
\parshape 20 0cm 12.8cm 0cm 12.8cm 6.0cm 6.8cm
6.0cm 6.8cm 6.0cm 6.8cm 6.0cm 6.8cm 6.0cm 6.8cm 6.0cm 6.8cm 6.0cm 6.8cm
6.0cm 6.8cm 6.0cm 6.8cm 6.0cm 6.8cm 6.0cm 6.8cm 6.0cm 6.8cm 6.0cm 6.8cm 
6.0cm 6.8cm 6.0cm 6.8cm 6.0cm 6.8cm 0cm 12.8cm 0cm 12.8cm
    Suppose that $\dim M\geqslant 3$. Let's consider the Newtonian
dynamical system with force field $\bold F$ on $M$. Suppose that
this system satisfies additional normality condition. Let's choose
some fixed point $p_0$ in $M$ and some nonzero vector $\bold v_0$ at
this point. \vadjust{\vskip -0.1cm\special{em:graph edis-05b.gif}
\vskip 0.1cm}Denote by $\alpha$ the hyperplane perpendicular to
vector $\bold v_0$, and denote by $\bold P$ orthogonal projector
to this hyperplane (see Fig\. 9.1). Then $\bold P$ is a value of
projector field $\bold P$ from \S\,5 at the point $q=(p_0,\bold v_0)$
of tangent bundle $TM$. One can draw infinitely many surfaces passing
through the point $p_0$ and being perpendicular to $\bold v_0$ at
this point. Let's choose one of them belonging to transformation
class and denote it by $S$. Then choose local coordinates $u^1,\,\ldots,
\,u^{n-1}$ in the neighborhood of the point $p_0$ and use this point
to write normalization condition \thetag{9.1}. According to the results
of \S\,7 for $n\geqslant 3$ weak the equalities \thetag{9.2} lead to
the system of Pfaff equations \thetag{7.16} for the function $\nu(p)$,
while arbitrariness of $\nu_0=|\bold v_0|$ in thetag{9.1} means that
these Pfaff equations are compatible in the sense of definition~8.1.
Hence the relationships \thetag{8.20} are fulfilled. In local
coordinates they are written as \thetag{8.21}.  We shall consider
them at the point $p_0$ on $S$.
\proclaim{Theorem 9.2} Let $S$ be some hypersurface in Riemannian manifold
$M$ and let $\bold P$ be operator of orthogonal projection onto the
hyperplane tangent to $S$ at the point $p_0$ on $S$. Then for the matrix
components of operator $\bold P$ in local coordinates $x^1,\,\ldots,\,
x^n$ on $M$ there exists the formula 
$$
P^i_\varepsilon=\sum^n_{s=1}\sum^{n-1}_{a=1}\sum^{n-1}_{k=1}
g_{\varepsilon s}\,\tau^s_a\,\rho^{ak}\,\tau^i_k,\hskip -2em
\tag9.6
$$
where $\tau^s_a$ and $\tau^i_k$ are components of coordinate tangent vectors
$\boldsymbol\tau_a$ and $\boldsymbol\tau_k$ that are defined by some local
coordinates $u^1,\,\ldots,\,u^{n-1}$ on $S$; $\rho^{ak}$ are components of
dual metric tensor for induced Riemannian metric $\bold\rho$ in coordinates
$u^1,\,\ldots,\,u^{n-1}$ on $S$; and finally, $g_{\varepsilon s}$ are
components of metric $\bold g$ in coordinates $x^1,\,\ldots,\,x^n$ on $M$.
\endproclaim
\demo{Proof} Let $\bold y$ be some arbitrary vector at the point $p_0$
on $M$. Its projection $\bold P\bold y$ can be expanded in the base of
vectors $\boldsymbol\tau_1,\,\ldots,\,\boldsymbol\tau_{n-1}$:
$$
\bold P\bold y=\sum^{n-1}_{k=1}\beta^k\,\boldsymbol\tau_k.\hskip -2em
\tag9.7
$$
In order to find numeric coefficients $\beta^k$ in \thetag{9.7} \pagebreak
we consider scalar products of the vector $\boldsymbol\tau_a$ with both
sides of the equality \thetag{9.7}. This yields
$$
(\boldsymbol\tau_a\,|\,\bold P\bold y)=\sum^{n-1}_{k=1}\beta^k\,
(\boldsymbol\tau_a\,|\,\boldsymbol\tau_k).\hskip -2em
\tag9.8
$$
In left hand side of \thetag{9.8} we obtain $(\boldsymbol\tau_a\,|\,
\bold P\bold y)=(\bold P\boldsymbol\tau_a\,|\,\bold y)=(\boldsymbol
\tau_a\,|\,\bold y)$, while in right hand side the quantities 
$(\boldsymbol\tau_a\,|\,\boldsymbol\tau_k)$ coincide with components
of metric tensor for induced Riemannian metric on $S$: $\rho_{ak}=
(\boldsymbol\tau_a\,|\,\boldsymbol\tau_k)$. Therefore we can bring
\thetag{9.8} to 
$$
(\boldsymbol\tau_a\,|\,\bold y)=\sum^{n-1}_{k=1}\beta^k\,\rho_{ak},
$$
and we can calculate $\beta^k$ with the use of components $\rho^{ak}$
for dual metric tensor:
$$
\beta^k=\sum^{n-1}_{a=1}(\boldsymbol\tau_a\,|\,\bold y)\,\rho^{ak}.
\hskip -2em
\tag9.9
$$
Let's substitute \thetag{9.9} into the formula \thetag{9.7} for
the vector $\bold P\bold y$. As a result we get
$$
\bold P\bold y=\sum^{n-1}_{a=1}\sum^{n-1}_{k=1}(\boldsymbol\tau_a\,|\,
\bold y)\,\rho^{ak}\,\boldsymbol\tau_k.\hskip -3em
\tag9.10
$$
Now we are only to write vectorial equality \thetag{9.10} in local
coordinates $x^1,\,\ldots,\,x^n$ on the manifold $M$. This yields
the following relationship:
$$
\sum^n_{\varepsilon=1}P^i_\varepsilon\,y^\varepsilon=\sum^n_{\varepsilon=1}
\sum^n_{s=1}\sum^{n-1}_{a=1}\sum^{n-1}_{k=1}g_{\varepsilon s}\,\tau^s_a\,
y^\varepsilon \,\,\rho^{ak}\,\tau^i_k.
$$
Components of vector $\bold y$ in both sides of the obtained equality
are arbitrary real numbers. Therefore it is reduced to the equality
\thetag{9.6} that was to be proved.\qed\enddemo
    We shall use theorem~9.2 to make further changes in \thetag{8.21}.
In order to do it let's take two copies of the equality \thetag{9.6},
redesignating indices in the second copy:
$$
\align
P^i_\varepsilon &=\sum^n_{s=1}\sum^{n-1}_{a=1}\sum^{n-1}_{k=1}
g_{\varepsilon s}\,\tau^s_a\,\rho^{ak}\,\tau^i_k,\\
P^j_\sigma &=\sum^n_{w=1}\sum^{n-1}_{b=1}\sum^{n-1}_{r=1}
g_{\sigma w}\,\tau^w_b\,\rho^{br}\,\tau^j_r.
\endalign
$$
By multiplying these two equalities, we get the relationship
$$
P^i_\varepsilon\,P^j_\sigma=\sum^n_{s=1}\sum^{n-1}_{a=1}\sum^{n-1}_{k=1}
\sum^n_{w=1}\sum^{n-1}_{b=1}\sum^{n-1}_{r=1}g_{\varepsilon s}\,\tau^s_a\,
\rho^{ak}\,g_{\sigma w}\,\tau^w_b\,\rho^{br}\,(\tau^i_k\,\tau^j_r).
$$
Here we intentionally enclosed into brackets last two terms in right hand
side of this relationships. Exactly the same terms $\tau^i_k\,\tau^j_r$
are in left hand side of \thetag{8.21}. Let's multiply \thetag{8.21}
by $g_{\varepsilon s}\,\tau^s_a\,\rho^{ak}\,g_{\sigma w}\,\tau^w_b\,
\rho^{br}$ and let's sum the resulting equality in $s$, $a$, $k$,
$w$, $b$, and $r$. Then left hand side of \thetag{8.21} takes the form
$$
\sum^n_{i=1}\sum^n_{j=1}P^i_\varepsilon\,P^j_\sigma\left(\,\shave{
\sum^n_{m=1}}n^m\,\frac{F_i\,\tilde\nabla_mF_j-F_j\,\tilde\nabla_mF_i}{\nu}+
\nabla_jF_i-\nabla_iF_j\right).
$$
Now let's introduce the following notations
$$
H^j_\varepsilon=\sum^n_{s=1}\sum^{n-1}_{a=1}\sum^{n-1}_{k=1}
\sum^{n-1}_{m=1}g_{\varepsilon s}\,\tau^s_a\,\rho^{ak}\,b^m_k\,
\tau^j_m.\hskip -3em
\tag9.11
$$
Then in right hand side of transformed equality \thetag{8.21} we can
detect the quantities $H^j_\varepsilon$ and $H^j_\sigma$ that were
formed as a result of the above transformations. I terms of notations
\thetag{9.11} right hand side of \thetag{8.21} is written as
$$
\sum^n_{i=1}\sum^n_{j=1}\nu\,(P^i_\varepsilon\,H^j_\sigma-
P^i_\sigma\,H^j_\varepsilon)\,\tilde\nabla_jF_i.
$$
Now, equating transformed expressions for left and right hand sides
of \thetag{8.21}, we complete first step in transforming this equality:
$$
\gathered
\sum^n_{i=1}\sum^n_{j=1}P^i_\varepsilon\,P^j_\sigma\left(\,\shave{
\sum^n_{m=1}}n^m\,\frac{F_i\,\tilde\nabla_mF_j-F_j\,\tilde\nabla_mF_i}
{\nu^2}\right.\,+\\
\vspace{2ex}
+\,\left.\frac{\nabla_jF_i-\nabla_iF_j}{\nu}\vphantom{\shave{\sum^n_{m=1}}
n^m\,\frac{F_i\,\tilde\nabla_mF_j-F_j\,\tilde\nabla_mF_i}{\nu}}\right)=
\sum^n_{i=1}\sum^n_{j=1}(P^i_\varepsilon\,H^j_\sigma-
P^i_\sigma\,H^j_\varepsilon)\,\tilde\nabla_jF_i.
\endgathered\hskip -2em
\tag9.12
$$\par
     Let's postpone for a while further transformation of \thetag{9.12}.
Instead, let's study the quantities $H^j_\varepsilon$ introduced by
\thetag{9.11}. They have transparent geometric interpretation. In
\thetag{9.11} we find the quantities $b^m_k$ that were derived from
$b_{kr}$ by index raising procedure (see formula \thetag{7.13}). They
are the components of inner tensor field $\bold B$ of type $(1,1)$
in $S$. The value of such field at the point $p_0\in S$ can be
interpreted as linear operator in the tangent space $T_{p_0}(S)$ to $S$.
Tangent space $T_{p_0}(S)$ is naturally embedded into the tangent space
$T_{p_0}(M)$ as $(n-1)$-dimensional hyperplane $\alpha$ (see Fig\. 9.1).
We have orthogonal projector $\bold P$ onto this hyperplane. Therefore
composition $\bold H=\bold B\compos\bold P$ can be interpreted as
linear operator in $T_{p_0}(M)$. Let's calculate the matrix of this
operator in local coordinates $x^1,\,\ldots,\,x^n$. Suppose that
$\bold y$ is some arbitrary tangent vector at the point $p_0$ on $M$.
Then
$$
\bold H(\bold y)=\bold B\compos\bold P(\bold y)=\bold B(\bold P\bold y).
\hskip -3em
\tag9.13
$$
In order to calculate projection $\bold P\bold y$ in \thetag{9.13} we
use formula \thetag{9.10}. The result of applying $\bold B$ to base
vectors $\boldsymbol\tau_1,\,\ldots,\,\boldsymbol\tau_{n-1}$ in hyperplane
$\alpha$ is defined by matrix of $\bold B$:
$$
\bold B(\boldsymbol\tau_k)=\sum^{n-1}_{m=1}b^m_k\,\boldsymbol\tau_m.
\hskip -3em
\tag9.14
$$
Combining \thetag{9.13}, \thetag{9.14}, and \thetag{9.10}, for the
vector $\bold H(\bold y)$ in \thetag{9.13} we get 
$$
\bold H(\bold y)=\sum^{n-1}_{m=1}\sum^{n-1}_{a=1}\sum^{n-1}_{k=1}
(\boldsymbol\tau_a\,|\,\bold y)\,\rho^{ak}\,b^m_k\,\boldsymbol\tau_m.
\hskip -3em
\tag9.15
$$
Now let's bring vectorial equality \thetag{9.15} to coordinate form:
$$
\sum^n_{\varepsilon=1}H^j_\varepsilon\,y^\varepsilon=
\sum^n_{\varepsilon=1}\sum^n_{s=1}\sum^{n-1}_{m=1}\sum^{n-1}_{a=1}
\sum^{n-1}_{k=1}g_{\varepsilon s}\,\tau^s_a\,y^\varepsilon\,
\rho^{ak}\,b^m_k\,\tau^j_m.
$$
Components of vector $\bold y$ are arbitrary real numbers. Taking into
account this fact, we find that components of operator $\bold H$ in local
coordinates $x^1,\,\ldots,\,x^n$ are given exactly by formula \thetag{9.11}.
\proclaim{Theorem 9.3} Operator $\bold H=\bold B\compos\bold P$ in
tangent space $T_{p_0}(M)$ defined by second fundamental form of
hypersurface $S$ according to the formula \thetag{9.11} for its
components is symmetric in Riemannian metric $\bold g$ of the
manifold $M$, i\.~e\.
$$
(\bold x\,|\,\bold H\bold y)=(\bold H\bold x\,|\,\bold y),
\hskip -3em
\tag9.16
$$
where $\bold x$ and $\bold y$ are arbitrary vectors from tangent space
$T_{p_0}(M)$ at the point $p_0\in S$.
\endproclaim
\demo{Proof} In order to prove the relationship \thetag{9.16} we use
the above formula \thetag{9.15} for the vector $\bold H\bold y=
\bold H(\bold y)$. This yields
$$
(\bold x\,|\,\bold H\bold y)=\sum^{n-1}_{m=1}\sum^{n-1}_{a=1}
\sum^{n-1}_{k=1}(\boldsymbol\tau_a\,|\,\bold y)\,\rho^{ak}\,b^m_k\,
(\boldsymbol\tau_m\,|\,\bold x).
$$
Taking into account \thetag{7.13}, this relationship can be transformed
as follows:
$$
(\bold x\,|\,\bold H\bold y)=\sum^{n-1}_{m=1}\sum^{n-1}_{a=1}
\sum^{n-1}_{k=1}\sum^{n-1}_{r=1}(\boldsymbol\tau_a\,|\,\bold y)
\,\rho^{ak}\,b_{kr}\,\rho^{rm}\,(\boldsymbol\tau_m\,|\,\bold x).
\hskip -3em
\tag9.17
$$
Right hand side of the equality \thetag{9.17} is invariant under the
transposition of vectors $\bold x$ and $\bold y$. This is due to the
symmetry in components of second fundamental form of $S$ (see equality
\thetag{7.12} in \S\,7). Therefore we get the required relationship
\thetag{9.16}: \pagebreak $(\bold x\,|\,\bold H\bold y)=(\bold y\,|\,
\bold H\bold x)=(\bold H\bold x\,|\,\bold y)$.\qed\enddemo

    Symmetric in metric $\bold g$ operator $\bold H$ commutate with
operator $\bold P$ of orthogonal projection to hyperplane $\alpha=
T_{p_0}(S)$, i\.~e\. we have:
$$
\bold P\compos\bold H=\bold H\compos\bold P=\bold H.\hskip -3em
\tag9.18
$$
Let $\bold x\in T_{p_0}(M)$. Operator $\bold H$ maps vector $\bold x$
to the vector $\bold H\bold x\in \alpha$. Therefore, applying $\bold P$
to $\bold H\bold x$, we get the same vector $\bold H\bold x$, i\.~e\.
$\bold P(\bold H\bold x)=\bold H\bold x$. This means that $\bold P\compos
\bold H=\bold H$. On the other hand $\bold H\compos\bold P=\bold B\compos
\bold P^2$. If we recall that $\bold P^2=\bold P$, then we get $\bold H
\compos\bold P=\bold B\compos\bold P=\bold H$. This proves the equality
\thetag{9.18}.
\proclaim{Theorem 9.4 (\bf on the second fundamental form)} Let $\bold v_0$
be nonzero vector at some point $p_0$ on Riemannian manifold $M$, and let
$\bold P$ be the operator of orthogonal projection onto the hyperplane
$\alpha$ perpendicular to $\bold v_0$. Then for any linear operator $\bold H$
in $T_{p_0}(M)$ symmetric in metric $\bold g$ of Riemannian manifold $M$
and satisfying the relationships \thetag{9.18} one can find hypersurface
$S$ passing through the point $p_0$ and perpendicular to $\bold v_0$ such
that matrix elements of the operator $\bold H$ are determined by second
fundamental form of $S$ according to the formula \thetag{9.11}.
\endproclaim
    We shall prove this theorem in \S\,10 (see below). Now we use it
to analyze the equations \thetag{9.12}. Identically zero operator
$\bold H=0$ satisfy the hypothesis of theorem~9.4. Therefore its 
matrix $H^i_j=0$ can be present in right hand side of \thetag{9.12}.
Substituting $H^j_\sigma=0$ and $H^j_\varepsilon=0$ into \thetag
{9.12} we obtain
$$
\sum^n_{i=1}\sum^n_{j=1}P^i_\varepsilon\,P^j_\sigma\left(
\nabla_jF_i-\nabla_iF_j+\shave{\sum^n_{m=1}}n^m\,\frac{F_i\,
\tilde\nabla_mF_j-F_j\,\tilde\nabla_mF_i}{\nu}\right)=0.
$$
Left  hand  side  of  \thetag{9.12}  doesn't  depend  on   second 
fundamental
form of $S$, hence it doesn't depend on $\bold H$. The above equality
means that left hand side of \thetag{9.12} is equal to zero. Therefore
right hand side of \thetag{9.12} is zero too:
$$
\sum^n_{i=1}\sum^n_{j=1}(P^i_\varepsilon\,H^j_\sigma-P^i_\sigma\,
H^j_\varepsilon)\,\tilde\nabla_jF_i=0.\hskip -3em
\tag9.19
$$
First of the above two relationships can be written in the form that
doesn't depend on $S$ at all. In order to do it let's recall that
normal vector to $S$ at the point $p_0$ is directed along the vector
of velocity, while function $\nu$ iz normalized at this point by
$\nu=\nu_0=|\bold v_0|$. Hence $\bold n$ can be replaced by $\bold N$,
and $\nu$ can be replaced by $v=|\bold v|$:
$$
\sum^n_{i=1}\sum^n_{j=1}P^i_\varepsilon\,P^j_\sigma\left(
\nabla_jF_i-\nabla_iF_j+\shave{\sum^n_{m=1}}N^m\,\frac{F_i\,
\tilde\nabla_mF_j-F_j\,\tilde\nabla_mF_i}{v}\right)=0.
$$
Arbitrariness in the choice of point $p_0$ (see definition~9.1) means
that this equality holds at all points $q=(p,\bold v)$ of tangent
bundle $TM$, where $\bold v\neq 0$. It is the differential equation
for the force field of dynamical system on $M$ that follows from
additional normality condition for it. For the sake of more symmetry
we shall write it as
$$
\gathered
\sum^n_{i=1}\sum^n_{j=1}P^i_\varepsilon\,P^j_\sigma\left(
\,\shave{\sum^n_{m=1}}N^m\,\frac{F_i\,\tilde\nabla_mF_j}{v}
-\nabla_iF_j\right)=\\
\vspace{1ex}
=\sum^n_{i=1}\sum^n_{j=1}P^i_\varepsilon\,P^j_\sigma\left(
\,\shave{\sum^n_{m=1}}N^m\,\frac{F_j\,\tilde\nabla_mF_i}{v}
-\nabla_jF_i\right).
\endgathered\hskip -3em
\tag9.20
$$\par
    The equation \thetag{9.19}, unlike to \thetag{9.20}, require
further transformation, since it contain the quantities $H^j_\sigma$
and $H^j_\varepsilon$, which depend on the choice of hypersurface $S$.
Let's consider a linear operator in $T_{p_0}(M)$ with the following
components:
$$
K^\sigma_\varepsilon=\sum^n_{m=1}\sum^n_{i=1}\sum^n_{j=1}
g^{\sigma m}\,P^j_m\,\tilde\nabla_jF_i\,P^i_\varepsilon.
\hskip -3em
\tag9.21
$$
Operator $\bold K$ with components \thetag{9.21} satisfies the
relationships
$$
\bold K\compos\bold P=\bold P\compos\bold K=\bold K.
\hskip -3em
\tag9.22
$$
The relationships \thetag{9.22} are analogous to the relationships
\thetag{9.18} for the operator $H$. They can be verified by direct
calculations, if we take into account that $\bold P^2=\bold P$ and
take into account the equalities
$$
\sum^n_{m=1}g^{\sigma m}\,P^j_m=\sum^n_{m=1}P^\sigma_m\,g^{mj}
\hskip -3em
\tag9.23
$$
that reflect the symmetry of projector $\bold P$ with respect to
Riemannian metric $\bold g$ on $M$. Further let's consider the
operator $\bold M=\bold H\compos\bold K$, and let's calculate its
components in local coordinates. For the components of $\bold M$
we obtain
$$
M^\gamma_\varepsilon=\sum^n_{\sigma=1}H^\gamma_\sigma\,
K^\sigma_\varepsilon=\sum^n_{\sigma=1}\sum^n_{i=1}\sum^n_{j=1}
g^{\gamma\sigma}\,H^j_\sigma\,\tilde\nabla_jF_i\,P^i_\varepsilon.
$$
In deriving this formula we used \thetag{9.23}, and used one of
the relationships \thetag{9.18}. By means of $\bold M$ we define
bilinear form
$$
\boldsymbol\Theta(\bold x,\bold y)=(\bold x\,|\,\bold M\bold y).
\hskip -3em
\tag9.24
$$
Form \thetag{9.24} is associated with tensor $\boldsymbol\Theta$
of type $(0,2)$. Here are its components:
$$
\theta_{\sigma\varepsilon}=\sum^n_{i=1}\sum^n_{j=1}H^j_\sigma\,
\tilde\nabla_jF_i\,P^i_\varepsilon.\hskip -3em
\tag9.25
$$
Let's compare formulas \thetag{9.25} and \thetag{9.19}. It's easy to
note that \thetag{9.19} is exactly the condition of symmetry for bilinear
form $\boldsymbol\Theta$, \pagebreak i\.~e\. $\theta_{\sigma\varepsilon}=
\theta_{\varepsilon\sigma}$ or $\boldsymbol\Theta(\bold x,\bold y)=
\boldsymbol\Theta(\bold y,\bold x)$. Hence $\bold M$ is a symmetric
operator in Riemannian metric $\bold g$ of the manifold $M$:
$$
(\bold x\,|\,\bold M\bold y)=(\bold M\bold x\,|\,\bold y),
\hskip -3em
\tag9.26
$$
Let's state this result as the following lemma.
\proclaim{Lemma 9.1} The equations \thetag{9.19} are equivalent to
the symmetry of operator $\bold M=\bold H\compos\bold K$ in the metric
$\bold g$ of $M$, operators $\bold H$ and $\bold K$ being defined by
their components \thetag{9.11} and \thetag{9.21}.
\endproclaim
     First let's take $\bold H=\bold P$ and let's apply theorem~9.4
to such operator $\bold H$. This is correct, since operator $\bold H=
\bold P$ is symmetric, and it satisfies the relationships \thetag{9.18}.
From theorem~9.4 we get that operator $\bold H=\bold P$ can be defined
by second fundamental form of some hypersurface $S$ passing through
the point $p_0$. For $\bold H=\bold P$ we obtain $\bold M=\bold K$.
Therefore due to lemma~9.1 we conclude that operator $\bold K$ with
components \thetag{9.21} is symmetric in Riemannian metric $\bold g$
of the manifold $M$.\par
     The equality $\bold H=\bold P$, which holds for some particular
choice of hypersurface $S$, is only some special case, when theorem~9.4
is applicable. In general, operator $\bold H$ is an arbitrary symmetric
operator satisfying the relationships \thetag{9.18} and such that
composition of two symmetric operators $\bold H$ and $\bold K$ is
symmetric operator $\bold M=\bold H\compos\bold K$. This is possible
if and only if operators $\bold H$ and $\bold K$ are commutating.
Indeed, from the equality \thetag{9.26} we derive
$$
(\bold x\,|\,\bold H(\bold K\bold y))=(\bold H(\bold K\bold x)\,|\,
\bold y).\hskip -3em
\tag9.27
$$
On the other hand symmetry of operators $\bold H$ and $\bold K$ yields
$$
(\bold x\,|\,\bold H(\bold K\bold y))=(\bold H\bold x\,|\,\bold K
\bold y)=(\bold K(\bold H\bold x)\,|\,\bold y).\hskip -3em
\tag9.28
$$
Comparing \thetag{9.27} and \thetag{9.28}, and taking into account
the arbitrariness of $\bold x$ and $\bold y$ in these relationships,
we obtain the operator equality
$$
\bold H\compos\bold K=\bold K\compos\bold H.\hskip -3em
\tag9.29
$$
This equality \thetag{9.29} means that $\bold H$ and $\bold K$ are
commutating. We state this result in form of the following lemma.
\proclaim{Lemma 9.2} Suppose that the equations \thetag{9.19} are
fulfilled for any hypersurface $S$ passing through the point $p_0$
and perpendicular to the vector $\bold v_0$ at this point. Then
operator $\bold K$ commutates with any symmetric operator $\bold H$
that satisfies the relationships \thetag{9.18}.
\endproclaim
    Let $\alpha$ be a hyperplane perpendicular to the vector
$\bold v_0$. Tangent space to the manifold $M$ at the point
$p_0$ is represented by a sum of two subspaces
$$
\pagebreak
T_{p_0}(M)=\alpha\oplus\langle\bold v_0\rangle.\hskip -3em
\tag9.30
$$
Here $\langle\bold v_0\rangle$ is a linear span of vector
$\bold v_0$. Each subspace in the expansion \thetag{9.30}
is invariant under the action of operators $\bold H$ and
$\bold K$. This is easily derived from \thetag{9.18} and
\thetag{9.22}. Moreover, the restrictions of $\bold H$ and
$\bold K$ to $\langle\bold v_0\rangle$ are zero. Let's
choose some orthonormal base $\bold E_1,\,\ldots,\,\bold E_{n-1}$
in hyperplane $\alpha$ and let's complete it by unitary vector
$\bold E_n$ directed along the vector $\bold v_0$. Operators
$\bold H$ and $\bold K$ in the base $\bold E_1,\,\ldots,\,\bold E_n$
are defined by the following symmetric matrices:
$$
\gather
H=\Vmatrix
H^1_1 & \hdots & H^1_{n-1} & 0\\
\vspace{1ex}
\vdots &\ddots & \vdots &\vdots\\
\vspace{2ex}
H^{n-1}_1 & \hdots & H^{n-1}_{n-1} & 0\\
\vspace{2ex}
0 &\hdots & 0 & 0
\endVmatrix,\hskip -3em
\tag9.31\\
\vspace{3ex}
K=\Vmatrix
K^1_1 & \hdots & K^1_{n-1} & 0\\
\vspace{1ex}
\vdots &\ddots & \vdots &\vdots\\
\vspace{2ex}
K^{n-1}_1 & \hdots & K^{n-1}_{n-1} & 0\\
\vspace{2ex}
0 &\hdots & 0 & 0
\endVmatrix.\hskip -3em
\tag9.32
\endgather
$$
Matrix of projector $\bold P$ in the base $\bold E_1,\,\ldots,\,
\bold E_n$ has similar structure:
$$
P=\Vmatrix
1 & \hdots & 0 & 0\\
\vspace{0.3ex}
\vdots &\ddots & \vdots &\vdots\\
\vspace{0.3ex}
0 & \hdots & 1 & 0\\
\vspace{0.3ex}
0 &\hdots & 0 & 0
\endVmatrix.\hskip -3em
\tag9.33
$$
Operator $\bold K$ is some fixed operator with matrix \thetag{9.32},
but $\bold H$, as it follows from lemma~9.2, is an arbitrary operator
with matrix of the form \thetag{9.31}. Using this arbitrariness, let's
choose the matrix of operator $\bold H$ to be the diagonal matrix with
distinct elements $H^1_1,\,\ldots,\,H^{n-1}_{n-1},\,H^n_n=0$:
$$
H=\Vmatrix
H^1_1 & \hdots & 0 & 0\\
\vspace{1ex}
\vdots &\ddots & \vdots &\vdots\\
\vspace{2ex}
0 & \hdots & H^{n-1}_{n-1} & 0\\
\vspace{2ex}
0 &\hdots & 0 & 0
\endVmatrix.\hskip -3em
\tag9.34
$$
Commutativity of operators $\bold K$ and $\bold H$ means that their
matrices \thetag{9.32} and \thetag{9.34} are also commutating. This
yields $K^i_j\,H^j_j=H^i_i\,K^i_j$. Diagonal elements in matrix
\thetag{9.34} are distinct: $H^i_i\neq H^j_j$ for $i\neq j$. Therefore
$K^i_j=0$ for $i\neq j$, i\.~e\. matrix of operator $\bold K$ is
diagonal in the base of vectors $\bold E_1,\,\ldots,\,\bold E_n$.
This reduces matrix \thetag{9.32} to the following form:
$$
K=\Vmatrix
K^1_1 & \hdots & 0 & 0\\
\vspace{1ex}
\vdots &\ddots & \vdots &\vdots\\
\vspace{2ex}
0 & \hdots & K^{n-1}_{n-1} & 0\\
\vspace{2ex}
0 &\hdots & 0 & 0
\endVmatrix.\hskip -3em
\tag9.35
$$
Now let's use again the arbitrariness of operator $\bold H$. In this
case we choose it so that its matrix \thetag{9.31} has the form:
$$
H=\Vmatrix
0 & \hdots & 1 &\hdots & 0\\
\vspace{0.3ex}
\vdots &\ddots & \vdots & &\vdots\\
\vspace{0.3ex}
1 & \hdots & 0 &\hdots & 0\\
\vspace{0.3ex}
\vdots & & \vdots &\ddots &\vdots\\
\vspace{0.3ex}
0 &\hdots & 0 &\hdots & 0
\endVmatrix.\hskip -3em
\tag9.36
$$
The only nonzero elements in matrix \thetag{9.36} are $H^1_i=1$ and
$H^i_1=1$, where $1<i<n$. Commutativity of $\bold K$ and $\bold H$
implies commutativity of their matrices \thetag{9.35} and \thetag{9.36}.
This yields $K^1_1\,H^1_i=H^1_i\,K^i_i$ and $K^i_i\,H^i_1=H^i_1\,K^1_1$. 
These relationships are equivalent to $K^i_i=K^1_1$. Due to arbitrariness
in $i$ we get
$$
K^1_1=\ldots=K^{n-1}_{n-1}=\lambda.\hskip -3em
\tag9.37
$$
The relationship \thetag{9.37} reduces matrix \thetag{9.35} to the
following form:
$$
K=\Vmatrix
\lambda & \hdots & 0 & 0\\
\vspace{0.3ex}
\vdots &\ddots & \vdots &\vdots\\
\vspace{0.3ex}
0 & \hdots & \lambda & 0\\
\vspace{0.3ex}
0 &\hdots & 0 & 0
\endVmatrix.
$$
Now it's obvious that such matrix commutates with any matrix of the form
\thetag{9.31}. Comparing it with \thetag{9.33}, we get $\bold K=\lambda
\cdot\bold P$. We state this result as a third lemma. 
\proclaim{Lemma 9.3} Operator $\bold K$ satisfying the relationships
\thetag{9.22} commutates with arbitrary symmetric operator $\bold H$
satisfying the relationships \thetag{9.18} if and only if it differs
from projector $\bold P$ only by some numeric multiple: $\bold K=\lambda
\cdot\bold P$.
\endproclaim
    Thus, by means of theorem~9.4 and by means of three lemmas~9.1, 9.2,
and 9.3 we reduced \thetag{9.19} to the operator equality
$$
\bold K=\lambda\cdot\bold P.\hskip -3em
\tag9.38
$$
Let's write \thetag{9.38} in local coordinates. In order to do it we
use formula \thetag{9.21} that determines components of operator
$\bold K$. As a result we get the following relationship:
$$
\sum^n_{m=1}\sum^n_{i=1}\sum^n_{j=1}g^{\sigma m}\,P^j_m\,
\tilde\nabla_jF_i\,P^i_\varepsilon=\lambda\,
P^\sigma_\varepsilon.\hskip -3em
\tag9.39
$$
In formula \thetag{9.39} one can lower index $\sigma$ and one can
raise index $\varepsilon$. If we take into account symmetry of operator
$\bold P$, we can write 
$$
\sum^n_{i=1}\sum^n_{j=1}P^j_\sigma\,\tilde\nabla_jF^i\,
P^\varepsilon_i=\lambda\,P^\varepsilon_\sigma.\hskip -3em
\tag9.40
$$
In the equations \thetag{9.40} we have numeric parameter $\lambda$.
Let's find it by calculating traces of operators in both sides of
the equality \thetag{9.38}:
$$
\tr\bold K=\lambda\cdot\tr\bold P.\hskip -3em
\tag9.41
$$
It's easy to calculate trace of projector: $\tr\bold P=n-1$. It is
nonzero, since we consider the case $n\geqslant 3$. Finding $\lambda$
from \thetag{9.41}, we can write \thetag{9.40} in the form that
has no indefinite parameters:
$$
\sum^n_{i=1}\sum^n_{j=1}P^j_\sigma\,\tilde\nabla_jF^i\,
P^\varepsilon_i=\sum^n_{i=1}\sum^n_{j=1}\sum^n_{m=1}
\frac{P^j_m\,\tilde\nabla_jF^i\,P^m_i}{n-1}\,P^\varepsilon_\sigma.
\hskip -3em
\tag9.42
$$\par
    Similar to \thetag{9.20}, the equations \thetag{9.42} are the
partial differential equations for the force field $\bold F$ of
Newtonian dynamical system. The above calculations can be reverted.
From \thetag{9.42} we could get \thetag{9.38}. Then from \thetag{9.38}
and \thetag{9.22} we could get the equality \thetag{9.29}, which
now is written as
$$
\bold H\compos\bold K=\bold K\compos\bold H=\lambda\cdot\bold H.
$$
Hence operator $\bold M=\bold H\compos\bold K=\lambda\cdot\bold H$ is
symmetric. Due to lemma~9.1 this is equivalent to \thetag{9.19}.
Combining \thetag{9.19} and \thetag{9.20}, we get the equality
\thetag{9.12}, which, in turn, is equivalent to \thetag{8.21}. The
equality \thetag{8.21} is the compatibility condition for the system
of Pfaff equations \thetag{7.16}. It provides solvability of these
equations for arbitrary choice of the point $p_0\in S$, and for
arbitrary choice of numeric parameter $\nu_0\neq 0$ in normalization
condition \thetag{9.1}. As a result we can formulate the main
theorem of this section, which was proved by the above considerations.
\proclaim{Theorem 9.5} Newtonian dynamical system on Riemannian
manifold of the dimension $n\geqslant 3$ satisfies additional
normality condition if and only if its force field satisfies the
equations \thetag{9.20} and \thetag{9.42} at all points $q=(p,\bold v)$
of tangent bundle $TM$, except for those, where $\bold v=0$.
\endproclaim
    The equations \thetag{9.20} and \thetag{9.42} are called
{\i additional normality conditions}. They were first derived
in \cite{38} (see also preprint \cite{34}) for the case
$M=\Bbb R^n$. \pagebreak The above derivation of the equations
\thetag{9.20} and \thetag{9.42} is based on the theorem~9.4,
which is not yet proved.
\head
\S\,10. Proof of the theorem on a second quadratic form.
\endhead
    In theorem~9.4 on the second fundamental form we consider a
point $p_0$ on Riemannian manifold $M$ and nonzero vector $\bold v_0$
at this point. Vector $\bold v_0$ defines a hyperplane $\alpha$ in
$T_{p_0}(M)$ perpendicular to $\bold v_0$ and an operator $\bold P$
of orthogonal projection onto the hyperplane $\alpha$. Then various
hypersurfaces passing through $p_0$ and tangent to $\alpha$ at that
point are considered. Basic tangent vectors $\boldsymbol\tau_1,\,
\ldots,\,\boldsymbol\tau_{n-1}$ at the point $p_0$ for any of such
hypersurfaces are in the hyperplane $\alpha$ and form the base in it.
Theorem~9.4 (which we need to prove) asserts that for any symmetric
linear operator $\bold H$ in $T_{p_0}(M)$ satisfying the relationships
\thetag{9.18} its matrix is defined by vectors $\boldsymbol\tau_1,\,
\ldots,\,\boldsymbol\tau_{n-1}$ and by second fundamental form of
some hypersurface $S$ according to the formula \thetag{9.11}. Vector
$\bold v_0$ and hyperplane $\alpha$ define the expansion of tangent
space $T_{p_0}(M)$ into a direct sum of two subspaces:
$$
T_{p_0}(M)=\alpha\oplus\langle\bold v_0\rangle.\hskip -3em
\tag10.1
$$
WE have already considered this expansion (see \thetag{9.30}).
Let $\bold H$ be some arbitrary symmetric in metric $\bold g$
on $M$ linear operator $\bold H\!:T_{p_0}(M)\to T_{p_0}(M)$
satisfying the relationships \thetag{9.18}. From $\bold H=\bold P
\compos\bold H$ we find that subspace $\alpha$ is invariant under
the action of $\bold H$. From $\bold H=\bold H\compos\bold P$, in
turn, we get that subspace $\langle\bold v_0\rangle$ in \thetag{10.1}
is also invariant under the action of $\bold H$, the restriction
of $\bold H$ to $\langle\bold v_0\rangle$ being zero. Let's denote
by $\bold B$ the restriction of $\bold H$ to hyperplane $\alpha$:
$$
\bold B=\bold H\,\hbox{\vrule height 8pt depth 6pt
width 0.5pt}_{\,\alpha}.
$$
Such restriction is a symmetric linear operator $\bold B\!:\alpha\to\alpha$
such that $\bold H=\bold B\compos\bold P$. And conversely, for any
symmetric linear operator $\bold B\!:\alpha\to\alpha$ the composition
$\bold H=\bold B\compos\bold P$ is symmetric operator $\bold H\!:
T_{p_0}(M)\to T_{p_0}(M)$ satisfying the relationships \thetag{9.18}.
\par
    Let $\boldsymbol\tau_1,\,\ldots,\,\boldsymbol\tau_{n-1}$ be some
arbitrary base in hyperplane $\alpha$. Consider a matrix with the
following components
$$
b_{ij}=(\boldsymbol\tau_i\,|\,\bold B\boldsymbol\tau_j).
\hskip -3em
\tag10.2
$$
Matrix \thetag{10.2} is symmetric, i\.~e\. $b_{ij}=b_{ji}$, this is
due to the symmetry of operator $\bold B$. Matrix with components
\thetag{10.2} can be obtained from the matrix of operator $\bold B$
in the base $\boldsymbol\tau_1,\,\ldots,\,\boldsymbol\tau_{n-1}$
by means of lowering index procedure:
$$
b_{ij}=\sum^{n-1}_{m=1}\rho_{im}\,b^m_j.
\hskip -3em
\tag10.3
$$
Here $\rho_{im}$ \pagebreak are components of metric tensor for the metric
$\bold g$ restricted to the hyperplane $\alpha$. Formula \thetag{10.3}
for $b_{ij}$ can be inverted as follows:
$$
b^m_k=\sum^{n-1}_{i=1}\rho^{mi}\,b_{ik}.
\hskip -3em
\tag10.4
$$
Operator $\bold B$ is reconstructed by matrix \thetag{10.2} according
to the formula
$$
\bold B\bold y=\sum^{n-1}_{a=1}\sum^{n-1}_{k=1}\sum^{n-1}_{m=1}
(\boldsymbol\tau_a\,|\,\bold y)\,\rho^{ak}\,b^m_k\,\boldsymbol
\tau_m.\hskip -3em
\tag10.5
$$
Formula \thetag{10.5} and formula \thetag{10.4} inverting \thetag{10.3}
mean that defining symmetric operator $\bold B$ in hyperplane $\alpha$
is equivalent to fixing some base $\boldsymbol\tau_1,\,\ldots,\,
\boldsymbol\tau_{n-1}$ in $\alpha$ and choosing some symmetric
matrix $b_{ij}$. Substituting \thetag{10.5} into the formula $\bold H=
\bold B\compos\bold P$, for components of operator $\bold H$ we get
the expression, which exactly the same as in \thetag{9.11}. Therefore
theorem~9.4 is the consequence of the following more simple
auxiliary theorem.
\proclaim{Theorem 10.1} Suppose that in tangent space $T_{p_0}(M)$
at some point $p_0$ on Riemannian manifold of the dimension $n\geqslant
3$ we choose some nonzero vector $\bold v_0$, mark hyperplane
$\alpha$ perpendicular to $\bold v_0$, mark some base $\boldsymbol
\tau_1,\,\ldots,\,\boldsymbol\tau_{n-1}$ in $\alpha$, and choose some
symmetric $(n-1)\times (n-1)$ matrix $b$. Then there exists some
hypersurface $S$ passing through the point $p_0$ tangent to $\alpha$
and there exist some local coordinates $u^1,\,\ldots,\,u^{n-1}$ on $S$
such that vectors $\boldsymbol\tau_1,\,\ldots,\,\boldsymbol\tau_{n-1}$
are coordinate tangent vectors to $S$ at the point $p_0$, while $b$
coincides with the matrix of second fundamental form for $S$ at this
point.
\endproclaim
    Before proving the theorem let's consider the following lemma that
states the fact, which is well-known in geometry.
\proclaim{Lemma 10.1} For any point $p_0$ on Riemannian manifold $M$
there exist local coordinates $x^1,\,\ldots,\,x^n$ such that component
of metric connection $\Gamma^k_{ij}$ vanishes at the point $p_0$ in
these local coordinates.
\endproclaim
\demo{Proof} Suppose that $\tx^1,\,\ldots,\,\tx^n$ are some arbitrary
local coordinates in some neighborhood of the point $p_0$. Without loss
of generality we can assume that
$$
\tx^1(p_0)=\ldots=\tx^n(p_0)=0.
$$
Let $\tilde\Gamma^k_{ij}(p_0)$ be components of metric connection at the
point $p_0$ in these local coordinates. Let's define new local coordinates
$x^1,\,\ldots,\,x^n$ by the following formula
$$
\pagebreak
x^k=\tx^k+\frac{1}{2}\sum^n_{i=1}\sum^n_{j=1}\tilde\Gamma^k_{ij}(p_0)
\,\tx^i\,\tx^j.\hskip -3em
\tag10.6
$$
For the components of transition matrix $S$ at the point $p_0$ from
\thetag{10.6} we get
$$
S^k_i(p_0)=\frac{\partial x^k}{\partial\tx^i}\,\hbox{\vrule height 12pt
depth 9pt width 0.5pt}_{\,p_0}=\delta^k_i.
$$
This means that $S$ is unitary matrix at the point $p_0$. Matrix $T=S^{-1}$
then is also unitary matrix at the point $p_0$. In order to calculate
components of metric connection $\Gamma^k_{ij}$ in newly defined local
coordinates let's use formula \thetag{4.11} from Chapter
\uppercase\expandafter{\romannumeral 3}: 
$$
\tilde\Gamma^k_{ij}=\sum^n_{m=1}\sum^n_{a=1}\sum^n_{c=1} T^k_m\,S^a_i
\,S^c_j\,\Gamma^m_{ac}+\sum^n_{m=1} T^k_m\,\frac{\partial
S^m_i}{\partial\tx^j}.\hskip -3em
\tag10.7
$$
If we take into account that $S$ and $T$ are unitary matrices at the
point $p_0$, then for components of metric connection $\Gamma^k_{ij}(p_0)$
and $\tilde\Gamma^k_{ij}(p_0)$ from \thetag{10.7} we get
$$
\tilde\Gamma^k_{ij}(p_0)=\Gamma^k_{ij}(p_0)+\frac{\partial^2 x^k}
{\partial\tx^i\,\partial\tx^j}\,\hbox{\vrule height 12pt
depth 9pt width 0.5pt}_{\,p_0}.\hskip -3em
\tag10.8
$$
Differentiating \thetag{10.6}, we calculate second derivative in
\thetag{10.8}:
$$
\frac{\partial^2 x^k}{\partial\tx^i\,\partial\tx^j}\,
\hbox{\vrule height 12pt depth 9pt width 0.5pt}_{\,p_0}
=\tilde\Gamma^k_{ij}(p_0).\hskip -3em
\tag10.9
$$
Now from \thetag{10.8} and \thetag{10.9} we see that components of
metric connection $\Gamma^k_{ij}$ for newly constructed local
coordinates \thetag{10.6} are zero at the point $p_0$, as it was
stated in lemma. Lemma is proved.\qed\enddemo
{\bf Remark}. Suppose that at the point $p_0$ we have $n$ vectors
$\boldsymbol\tau_1,\,\ldots,\,\boldsymbol\tau_n$ forming the base
in tangent space $T_{p_0}(M)$. Local coordinates $x^1,\,\ldots,\,x^n$
with $\Gamma^k_{ij}(p_0)=0$ can be chosen so that corresponding
coordinate tangent vectors
$$
\bold E_1=\frac{\partial}{\partial x^1},\ \ldots,\ \bold E_n=
\frac{\partial}{\partial x^n}\hskip -4em
\tag10.10
$$
at the point $p_0$ will coincide with vectors $\boldsymbol\tau_1,\,
\ldots,\,\boldsymbol\tau_n$. Indeed, the condition of vanishing
$\Gamma^k_{ij}$ at the point $p_0$ is invariant under the linear
change of local coordinates
$$
\xalignat 2
&x^k=\sum^n_{i=1}S^k_i\,\tx^i,
&&\tx^i=\sum^n_{i=1}T^i_k\,x^k,
\endxalignat
$$
where $S$ and $T=S^{-1}$ are some constant matrices. While vectors
\thetag{10.10} under such linear change of local coordinates are
transformed as follows:
$$
\xalignat 2
&\tilde\bold E_i=\sum^n_{i=1}S^k_i\,\bold E_k,
&&\bold E_k=\sum^n_{i=1}T^i_k\,\tilde\bold E_i.
\endxalignat
\pagebreak
$$
Therefore we can bring vectors $\bold E_1,\,\ldots,\,\bold E_n$ to
the coincidence with $\boldsymbol\tau_1,\,\ldots,\,\boldsymbol\tau_n$,
keeping condition $\Gamma^k_{ij}(p_0)=0$ fulfilled.
\demo{Proof of theorem 10.1} In theorem~10.1 we consider the point
$p_0$ with nonzero vector $\bold v_0$ at this point and the base 
$\boldsymbol\tau_1,\,\ldots,\,\boldsymbol\tau_{n-1}$ in hyperplane
$\alpha\subset T_{p_0}(M)$ perpendicular to the vector $\bold v_0$.
By denoting
$$
\boldsymbol\tau_n=\frac{\bold v_0}{|\bold v_0|}
$$
we complete base $\boldsymbol\tau_1,\,\ldots,\,\boldsymbol\tau_{n-1}$
up to the base $\boldsymbol\tau_1,\,\ldots,\,\boldsymbol\tau_n$ in
$T_{p_0}(M)$. Let's apply lemma~10.1 and choose local coordinates
$x^1,\,\ldots,\,x^n$ in the neighborhood of $p_0$ on $M$ such that
$\Gamma^k_{ij}(p_0)=0$ and such that
$$
\bold E_1=\boldsymbol\tau_1,\ \ldots,\ \bold E_{n-1}=\boldsymbol
\tau_{n-1},\qquad\bold E_n=\boldsymbol\tau_n=\frac{\bold v_0}
{|\bold v_0|}.
$$
Moreover, without loss of generality we can assume that
$x^1(p_0)=\ldots=x^n(p_0)=0$. Now let's define hypersurface $S$
by the equation
$$
x^n=\frac{1}{2}\sum^{n-1}_{i=1}\sum^{n-1}_{j=1}b_{ij}\,x^i\,x^j
\hskip -4em
\tag10.11
$$
in local coordinates chosen above. Here $b_{ij}$ are components
of constant symmetric matrix $b$ from the statement of theorem~10.1.
We can write \thetag{10.11} in parametric form by introducing local
coordinates $u^1,\,\ldots,\,u^{n-1}$ on $S$:
$$
\align
x^k&=u^k\text{ \ for \ }k=1,\,\ldots,\,n-1,\hskip -4em\\
\vspace{-1.4ex}
&\tag10.12\\
\vspace{-1.4ex}
x^n&=\frac{1}{2}\sum^{n-1}_{i=1}\sum^{n-1}_{j=1}b_{ij}\,u^i\,u^j.
\endalign
$$
Let's differentiate parametric equations \thetag{10.12} with respect to
the variable $u^i$:
$$
\frac{\partial x^k}{\partial u^i}=1\text{ \ for \ }i=k,\quad
\frac{\partial x^k}{\partial u^i}=0\text{ \ for \ }i\neq k<n.
\hskip -4em
\tag10.13
$$
By differentiating $x^n$ with respect to $u^i$ we get the formula
$$
\frac{\partial x^n}{\partial u^i}=\sum^{n-1}_{j=1}b_{ij}\,u^j.
\hskip -4em
\tag10.14
$$
Hypersurface $S$ defined by the equations \thetag{10.12} in parametric
form passes through the point $p_0$; this point has zero local coordinates
on $S$: $u^1(p_0)=\ldots=u^{n-1}(p_0)=0$. Partial derivatives \thetag{10.13}
and \thetag{10.14} form components of $i$-th coordinate tangent vector
$\boldsymbol\tau_i$ of $S$ in local coordinates $x^1,\,\ldots,\,x^n$ on
$M$:
$$
\xalignat 2
&\quad\tau^k_i=
\cases
1 &\text{for \ }i=k,\\
0 &\text{for \ }i\neq k<n,
\endcases
&&\tau^n_i=\sum^{n-1}_{j=1}b_{ij}\,u^j.\hskip -2em
\tag10.15
\endxalignat
$$
Substituting coordinates $u^1(p_0)=\ldots=u^{n-1}(p_0)$ of the point $p_0$
into these formulas, we get $\tau^n_i=0$. This means that $S$ is tangent
to $\alpha$, while coordinate tangent vectors to $S$ at the point $p_0$
coincide with the base $\bold E_1=\boldsymbol\tau_1,\,\ldots,\,\bold
E_{n-1}=\boldsymbol\tau_{n-1}$ in hyperplane $\alpha$. This proves first
proposal in theorem~10.1.\par
    In order to prove second proposal in theorem~10.1 let's calculate
the matrix of second fundamental form of hypersurface $S$ at the point
$p_0$, applying derivational formula of Weingarten \thetag{7.9} for this
purpose:
$$
\nabla_{u^j}\boldsymbol\tau_i=\sum^{n-1}_{m=1}\theta^m_{ij}\,
\boldsymbol\tau_m+\beta_{ij}\,\bold n.\hskip -4em
\tag10.16
$$
Let's calculate components of vector $\nabla_{u^j}\boldsymbol\tau_i$
in left hand side of derivational formula \thetag{10.16} according
to formula \thetag{5.8} from Chapter
\uppercase\expandafter{\romannumeral 4}:
$$
\nabla_{u^j}\tau^k_i=\frac{\partial\tau^k_i}{\partial u^j}+
\sum^n_{m=1}\sum^n_{r=1}\tau^m_j\,\Gamma^k_{mr}\,\tau^r_i.\hskip -4em
\tag10.17
$$
At the point $p=p_0$ we have $\Gamma^k_{mr}(p_0)=0$ due to special choice
of local coordinates $x^1,\,\ldots,\,x^n$ on $M$. For $\tau^k_i$ we apply
formula \thetag{10.15}. Then
$$
\xalignat 2
&\nabla_{u^j}\tau^k_i(p_0)=0\text{ \ for \ }k<n,
&&\nabla_{u^j}\tau^n_i(p_0)=b_{ij}.
\endxalignat
$$
Let's write these relationships in vectorial form, and let's calculate
components of unitary normal vector $\bold n$ to $S$ at the point $p_0$:
$$
\xalignat 2
&\qquad\nabla_{u^j}\boldsymbol\tau_i(p_0)=\Vmatrix 0\\ \vdots\\ 0\\
b_{ij}\endVmatrix,&&\bold n(p_0)=\Vmatrix 0\\ \vdots\\ 0\\ 1\endVmatrix.
\hskip -2em
\tag10.18
\endxalignat
$$
The relationships \thetag{10.15} for the quantities $\tau^i_j$,
written in vectorial form, yield
$$
\boldsymbol\tau_1(p_0)=\Vmatrix 1\\ \vdots\\ 0\\ 0\endVmatrix,\
\ldots,\
\boldsymbol\tau_{n-1}(p_0)=\Vmatrix 0\\ \vdots\\ 1\\ 0\endVmatrix.
\hskip -4em
\tag10.19
$$
Now, substituting \thetag{10.18} and \thetag{10.19} into vectorial
equality \thetag{10.16}, we find the values of parameters $\theta^m_{ij}$
$\beta_{ij}$ at the point $p_0$:
$$
\xalignat 2
&\theta^k_{ij}(p_0)=0,
&&\beta_{ij}=b_{ij}.
\endxalignat
$$
So, we see that matrix of second fundamental form $\beta$ for the
hypersurface $S$ defined by the equation \thetag{10.11} at the point
$p_0$ coincides with matrix $b$. Theorem~10.1 is now proved.\qed\enddemo
    Along with theorem~10.1 we have proved theorem~9.4 on the second
fundamental form. This follows from considerations preceding the
statement of theorem~10.1.
\head
\S\,11. Newtonian dynamical systems admitting the normal shift.
\endhead
\rightheadtext{\S\,11. Newtonian dynamical systems \dots}
     Having proved the theorem on second quadratic form, let's return
to the results of \S\,6 and \S\,9. Central points of these two sections
are definitions~6.1 and 9.1 that introduce weak and additional normality
conditions for Newtonian dynamical systems on Riemannian manifolds.
We managed to write each of these two conditions in form of systems of
partial differential equations for the force field $\bold F$ of the
dynamical system (see equations \thetag{6.16} and \thetag{6.17} for
weak normality condition, and equations \thetag{9.20} and \thetag{9.42}
for additional normality condition).
\definition{Definition 11.1} Weak normality condition from definition~6.1
and additional normality condition from definition~9.1, combined together,
form {\bf complete normality} condition for Newtonian dynamical system with
force field $\bold F$ on Riemannian manifold $M$.
\enddefinition
    The equations \thetag{6.16}, \thetag{6.17}, \thetag{9.20}, and
\thetag{9.42} form complete system of normality equations for the force
field of dynamical system $\bold F$. The whole set of these equations
is equivalent to complete normality condition from definition~11.1. 
Identically zero field $\bold F=0$ satisfies all these equations, therefore
geodesic flow is a simples example of dynamical system that satisfies
complete normality condition.\par
    Let's consider some Newtonian dynamical system with force field
$\bold F$ satisfying complete normality condition. Let $S$ be some
hypersurface in the manifold $M$. Let's mark some point $p_0\in S$ and
choose some local coordinates $u^1,\,\ldots,\,u^{n-1}$ in the
neighborhood of the point $p_0$ on $S$. Suppose $S'$ to be some smaller
piece of hypersurface $S$ containing the point $p_0$. Let's draw the
trajectories of dynamical system
$$
\xalignat 2
&\dot x^k=v^k,&&\nabla_tv^k=F^k,\hskip -3em
\tag11.1
\endxalignat
$$
passing through each point $p$ of $S'$. We define them by
the following initial data:
$$
\xalignat 2
&\qquad x^k(t)\,\hbox{\vrule height 8pt depth 8pt width 0.5pt}_{\,t=0}
=x^k(p),
&&v^k(t)\,\hbox{\vrule height 8pt depth 8pt width 0.5pt}_{\,t=0}=
\nu(p)\cdot n^k(p).\hskip -3em
\tag11.2
\endxalignat
$$
This arranges the shift $f_t\!:S'\to S'_t$ of $S'\subset S$ along the
trajectories of the dynamical system \thetag{11.1}. Piece $S'$ in $S$
can always be chosen belonging to the transformation class (see
definition~2.1). Due to additional normality condition being the part
of complete normality condition we can choose the function $\nu(p)$ on
$S'$ such that the shift $f_t\!:S'\to S'_t$ approximates normal shift
up to a first order derivatives. This means that all functions of
deviation $\varphi_1,\,\ldots,\,\varphi_{n-1}$ satisfy initial conditions
$$
\xalignat 2
&\qquad \varphi_k\,\hbox{\vrule height 8pt depth 8pt width 0.5pt}_{\,t=0}
=0,
&&\dot\varphi_k\,\hbox{\vrule height 8pt depth 8pt width 0.5pt}_{\,t=0}=0
\hskip -3em
\tag11.3
\endxalignat
$$
on initial hypersurface $S'$ (see conditions \thetag{7.5} above).
Function $\nu$ thereby can be chosen belonging to transformation
class (see definition~2.2) and normalized by 
$$
\nu(p_0)=\nu_0,\hskip -3em
\tag11.4
$$
where $\nu_0$ is an arbitrary nonzero real number.\par
   Weak normality condition, which also is a part of complete normality
condition, means that all functions of deviation $\varphi_1,\,\ldots,\,
\varphi_{n-1}$ satisfy linear homogeneous differential equations of
the second order in $t$ on trajectories of shift:
$$
\ddot\varphi_k-\Cal A(t)\,\dot\varphi_k-\Cal B(t)\,\varphi_k=0.
\hskip -3em
\tag11.5
$$
Combining differential equations \thetag{11.5} and initial data
\thetag{11.5}, we find that functions of deviation $\varphi_1,\,
\ldots,\,\varphi_{n-1}$ for the shift $f_t\!:S'\to S'_t$ are
identically zero. This, in turn, means that the shift we are
considering now is a normal shift in the sense of definition~2.3.
\definition{Definition 11.2} Newtonian dynamical system on Riemannian
manifold $M$ is called the system {\bf admitting the normal shift}, if
for any hypersurface $S$ in $M$, for any point $p_0\in S$, and for any
real number $\nu_0\neq 0$ there is a piece $S'$ of hypersurface $S$
belonging to transformation class and containing $p_0$, and there exists
a function $\nu(p)$ belonging to transformation class on $S'$ and
normalized by \thetag{11.4} such that the shift $f_t\!:S'\to S'_t$
defined by this function is a normal shift along the trajectories
of dynamical system.
\enddefinition
    Definition~11.2 states the central concept of this thesis, the
concept of Newtonian dynamical system {\bf admitting the normal shift}.
The condition stated i this definition was first introduced in \cite{39}.
There it was called {\bf strong normality} condition. In earlier papers
\cite{34}, \cite{35}, \cite{38}, and \cite{58} we used more simple version
of this condition, which did not include normalization \thetag{11.4}.
This condition was called {\bf normality condition}. It is more natural
from geometric point of view, but it is not in a good agreement with the
Pfaff equations that arises in our theory. Therefore presently, saying
dynamical system {\bf admitting the normal shift}, we imply the
definition~11.2.
\head
\S\,12. Equivalence of strong and complete normality conditions.
\endhead
\rightheadtext{\S\,12. Equivalence of strong and complete
normality \dots}
     Complete normality condition, which consists of weak normality
condition from definition~6.1 and additional normality condition
from definition~9.1, \pagebreak is sufficient for the strong normality
condition from definition~11.2 to be fulfilled. This follows from above
considerations preceding the statement of definition~11.2. Thus,
Newtonian dynamical system satisfying both weak and additional normality
conditions appears to be admitting the normal shift. It is applicable for
arranging normal shift of any preassigned hypersurface. This result can
be strengthened in form of the following theorem.
\proclaim{Theorem 12.1} {\bf Complete} and {\bf strong} normality
conditions for Newtonian dynamical systems on Riemannian manifolds
are {\bf equivalent} to each other.
\endproclaim
\demo{Proof} As we showed in \S\,11 above, complete normality condition
is {\bf sufficient} for the strong normality condition to be fulfilled.
Now we are to show that it is {\bf necessary} as well. Suppose that strong
normality condition for Newtonian dynamical system \thetag{11.1} is
fulfilled. This means that for any hypersurface $S$, for any point
$p_0\in S$, and for any real number $\nu_0\neq 0$ there is a piece $S'$
of hypersurface $S$ and there exists a function $\nu(p)$ on $S'$
normalized by the condition \thetag{11.4} and such that all functions
of deviations on shift trajectories are identically zero:
$$
\varphi_k(t)=0\text{ \ for all \ }k=1,\,\ldots,\,n-1.
\hskip -3em
\tag12.1
$$
This is the very condition that provides normality of shift $f_t\!:S'\to
S'$. From \thetag{12.2} we get that functions of deviations and their
first derivatives are zero at initial instant of time, i\.~e\. the
initial conditions \thetag{11.3} are fulfilled. Hence strong normality
condition implies additional normality condition to be fulfilled.\par
     Now we shall show that strong normality condition implies weak
normality condition to be fulfilled as well. Let's fix some trajectory
$p(t)$ of Newtonian dynamical system \thetag{11.1} and let's mark
some point $p_0=p(0)$ on it. Suppose that velocity vector $\bold v_0
=\bold v(0)$ at this point is nonzero: $\bold v_0\neq 0$. Then in tangent
space $T_{p_0}(M)$ we can take hyperplane $\alpha$ perpendicular to the
vector $v_0$. Let's consider various hyperplanes $S$ in $M$ passing
through the point $p_0$. For our preassigned trajectory $p(t)$ to be
one of the trajectories of normal shift for $S$ it should be
perpendicular to $S$ at the point $p_0$, i\.~e\. $S$ should be tangent
to hyperplane $\alpha$ at this point. Let's define unitary normal vector
at $p_0$ and a real number $\nu_0\neq 0$:
$$
\xalignat 2
&\bold n=\frac{\bold v_0}{|\bold v_0|},
&&\nu_0=|\bold v_0|.
\endxalignat
$$
This defines smooth field of unitary normal vectors on $S$ in some
neighborhood of $p_0$. Relying upon strong normality condition of
dynamical system, let's choose some piece $S'$ of hypersurface $S$
containing $p_0$, and choose some function $\nu(p)$ on $S'$ that
satisfies the condition \thetag{11.4} and initiates normal shift
$f_t\!:S'\to S'_t$. Let $u^1,\,\ldots,\,u^{n-1}$ be local coordinates
on $S'$. These local coordinates define variation vectors $\boldsymbol
\tau_1,\,\ldots,\,\boldsymbol\tau_{n-1}$ and deviation functions
$$
\pagebreak
\varphi_k=(\boldsymbol\tau_i\,|\,\bold v)\text{, \ \ where\ \ \ }
k=1,\,\ldots,\,n-1.\hskip -3em
\tag12.2
$$
Time derivatives of variation functions \thetag{12.2} are calculated
by formula \thetag{4.7}. Here this formula is written in the following
form:
$$
\dot\varphi_k=(\bold F\,|\,\boldsymbol\tau_k)+(\bold v\,|\,
\nabla_t\boldsymbol\tau_k).\hskip -3em
\tag12.3
$$
Covariant derivatives $\nabla_t\boldsymbol\tau_k$ in \thetag{12.3} at
the initial instant of time $t=0$ are calculated by formula \thetag{7.14}.
Using this formula, we take into account that function $\nu(p)$ on
initial hypersurface $S'$ satisfies the equation \thetag{7.16}. Hence
for covariant derivatives $\nabla_t\boldsymbol\tau_k$
in formula \thetag{12.3} we obtain
$$
\nabla_t\boldsymbol\tau_k\,\hbox{\vrule height 8pt depth 8pt width
0.5pt}_{\,t=0}=-\frac{(\bold F\,|\,\boldsymbol\tau_k)}{\nu}\cdot
\bold n-\nu\cdot\sum^{n-1}_{m=1}b^m_k\,\boldsymbol\tau_m.\hskip -3em
\tag12.4
$$
Normal vector $\bold n$ of hypersurface $S'$ is directed along the
vector of velocity; the function $\nu(p)$ determines the modulus of
initial velocity on $S'$. Therefore 
$$
\nabla_t\boldsymbol\tau_k\,\hbox{\vrule height 8pt depth 8pt width
0.5pt}_{\,t=0}=-\frac{(\bold F\,|\,\boldsymbol\tau_k)}{|\bold v|^2}
\cdot\bold v-|\bold v|\cdot\sum^{n-1}_{m=1}b^m_k\,\boldsymbol\tau_m.
\hskip -3em
\tag12.5
$$
Strong normality condition, that provides normality of shift $f_t\!:S'\to
S'_t$, implies identical in $t$ vanishing of all functions of
deviations. In particular this means 
$$
\xalignat 2
&\qquad \varphi_k\,\hbox{\vrule height 8pt depth 8pt width 0.5pt}_{\,t=0}
=0,
&&\dot\varphi_k\,\hbox{\vrule height 8pt depth 8pt width 0.5pt}_{\,t=0}=0.
\hskip -3em
\tag12.6
\endxalignat
$$
Vanishing of $\varphi_k$ and $\dot\varphi_k$ in \thetag{12.6} is
granted by orthogonality $\boldsymbol\tau_k\perp\bold v$ upon
substituting \thetag{12.5} into \thetag{12.3}. Let's consider
second derivatives 
$$
\ddot\varphi_k\,\hbox{\vrule height 8pt depth 8pt width 0.5pt}_{\,t=0}=0.
\hskip -3em
\tag12.7
$$
Their vanishing is also the consequence of $\varphi_k(t)=0$. We 
calculate second order derivatives in the left hand side of
\thetag{12.7} by formula \thetag{5.14}. As a result we get
$$
\boldsymbol\beta(\boldsymbol\tau_k)-\frac{\boldsymbol\alpha(\bold N)}
{|\bold v|}\,(\bold F\,|\,\boldsymbol\tau_k)-|\bold v|\sum^{n-1}_{m=1}
b^m_k\,\boldsymbol\alpha(\boldsymbol\tau_m)=0.
\hskip -3em
\tag12.8
$$
Here covectorial fields $\boldsymbol\alpha$ and $\boldsymbol\beta$ are
defined by their components \thetag{5.13}. When referred to the marked
point $p_0$, the quantities $b^m_k$ can be interpreted as components of
symmetric operator $\bold B$ in hyperplane $\alpha$ orthogonal to the
vector of velocity $\bold v=\bold v_0$ (see formula \thetag{9.14}).
Therefore equations \thetag{12.8} can be written as
$$
\boldsymbol\beta(\boldsymbol\tau_k)-\frac{\boldsymbol\alpha(\bold N)}
{|\bold v|}\,(\bold F\,|\,\boldsymbol\tau_k)-|\bold v|\,\boldsymbol
\alpha(\bold B\,\boldsymbol\tau_k)=0.
\hskip -3em
\tag12.9
$$
For further analysis of the equations \thetag{12.9} we use the
arbitrariness in the choice of hypersurface $S$ passing through
the point $p_0$ and perpendicular to vector $\bold v_0$ at this
point. Let's apply theorem~10.1. According to this theorem vector
$\boldsymbol\tau_k$ can be replaced by arbitrary vector $\boldsymbol
\tau$ from hyperplane $\alpha$, while operator $\bold B$ in \thetag{12.9}
can be understood as arbitrary symmetric operator $\bold B\!:\alpha
\to\alpha$. Hence vector $\bold B\,\boldsymbol\tau_k$ in \thetag{12.9}
can be replaced by some other arbitrary vector from hyperplane $\alpha$,
which do not depend on $\boldsymbol\tau$. Due to these reasons
equation \thetag{12.9} breaks into two parts
$$
\xalignat 2
&\boldsymbol\beta(\boldsymbol\tau)-\frac{\boldsymbol\alpha(\bold N)}
{|\bold v|}\,(\bold F\,|\,\boldsymbol\tau)=0,
&&\boldsymbol\alpha(\boldsymbol\tau)=0,
\endxalignat
$$
where $\boldsymbol\tau\in\alpha$ is arbitrary vector perpendicular to
$\bold v_0$ at the point $p_0$. We can avoid restriction $\boldsymbol\tau
\perp\bold v_0$, if we substitute $\boldsymbol\tau$ by $\bold P\boldsymbol
\tau$:
$$
\xalignat 2
&\boldsymbol\beta(\bold P\boldsymbol\tau)=\frac{\boldsymbol\alpha(\bold N)}
{|\bold v|}\,(\bold F\,|\,\bold P\boldsymbol\tau),
&&\boldsymbol\alpha(\bold P\boldsymbol\tau)=0.\hskip -4em
\tag12.10
\endxalignat
$$
Here $\boldsymbol\tau$ is arbitrary vector at the point $p_0$ on $M$.
Second equation \thetag{12.10} coincides with first equation \thetag{6.11},
since $\boldsymbol\tau$ in \thetag{12.10} and $\nabla_t\boldsymbol\tau$
in \thetag{6.11} stand for arbitrary vectors at marked point $p_0$ on
the trajectory of shift. First equation \thetag{12.10} coincides with
second equation \thetag{6.11}. In order to see it one should only denote
$$
\xalignat 2
&\quad\Cal A=\frac{\boldsymbol\alpha(\bold N)}{|\bold v|},
&&\Cal B=\frac{\boldsymbol\beta(\bold N)}{|\bold v|}-\Cal A\,
\frac{(\bold F\,|\,\bold N)}{|\bold v|}\hskip -4em
\tag12.11
\endxalignat
$$
(compare with \thetag{6.10} in \S\,6). Now let's use arbitrariness in
the choice of point $p_0$ and in the choice of vector $\bold v_0$ at
this point. This means that the equations \thetag{12.10} coinciding
with weak normality equations \thetag{6.17} and \thetag{6.17} are
fulfilled at all points $q=(p,\bold v)$  of tangent bundle $TM$,
where $\bold v\neq 0$.\par
     Now let's take some arbitrary function of deviation $\varphi$
on the trajectory $p(t)$ and form linear combination of its time
derivatives
$$
\ddot\varphi-\Cal A\,\dot\varphi-\Cal B\,\varphi,\hskip -4em
\tag12.12
$$
using parameters $\Cal A$ and $\Cal B$ from \thetag{12.11} as
coefficients in it. Here we do not associate vector of variation
$\boldsymbol\tau$ in $\varphi=(\bold v\,|\,\boldsymbol\tau)$ with
any hypersurface, and therefore we do not restrict it by the
condition $\bold v\perp\boldsymbol\tau$, as it was in the case
of normal shift. Here $\boldsymbol\tau=\boldsymbol\tau(t)$ is
arbitrary vector function that varies in $t$ according to the
differential equation \thetag{3.10}. From \thetag{3.10} we derive
the equalities \thetag{5.12} and \thetag{5.14} for $\dot\varphi$
and $\ddot\varphi$. Substituting them into \thetag{12.12} and taking
into account \thetag{12.10} and \thetag{12.11}, we find that
linear combination is identically zero: $\ddot\varphi-\Cal A\,\dot\varphi
-\Cal B\,\varphi=0$. This means that arbitrary function of deviation
satisfies homogeneous linear ordinary differential equation of the
second order with respect to parameter $t$ on trajectories of dynamical
system. According to definition~6.1 this is exactly the weak normality
condition for that dynamical system. Theorem~12.1 is proved.
\enddemo
\enddocument
\end